\documentclass{amsart}
\usepackage[latin1]{inputenc}
\usepackage{tikz}
\usetikzlibrary{arrows}
\usepackage{color} 
\usepackage{hyperref}
\hypersetup{colorlinks=true,linktoc=all,linkcolor=blue,citecolor=red,filecolor=pink,urlcolor=blue}
\usepackage[hyperref=true,style=numeric,backend=bibtex,maxnames=100,doi=false,isbn=false,url=false]{biblatex}
\bibliography{biblio}
\usepackage{csquotes}

\usepackage{amssymb,amsmath,latexsym}
\theoremstyle{plain}                    
\newtheorem{theorem}{Theorem}[section]
\newtheorem{lemma}[theorem]{Lemma}
\newtheorem{proposition}[theorem]{Proposition}
\newtheorem{corollary}[theorem]{Corollary}
\theoremstyle{definition}
\newtheorem{definition}[theorem]{Definition}
\newtheorem{example}[theorem]{Example}
\theoremstyle{remark}
\newtheorem{rmk}[theorem]{Remark}

\numberwithin{equation}{section}

\newcommand{\real}{\mathbb{R}}
\newcommand{\cmp}{\mathbb{C}}

\newcommand{\nat}{\mathbb{N}}

\newcommand{\hyp}{\mathbb{H}}

\newcommand{\cp}{\mathbb{C}\mathbb{P}^1}
\newcommand{\rp}{\mathbb{R}\mathbb{P}^1}
\newcommand{\pslc}{\mathrm{PSL}_2\cmp}

\newcommand{\pslr}{\mathrm{PSL}_2\real}

\newcommand{\qf}{Fuchsian }

\begin{document}
\title[BUBBLING $\cp$-STRUCTURES WITH QUASI-FUCHSIAN HOLONOMY]{BUBBLING COMPLEX PROJECTIVE STRUCTURES\\ WITH QUASI-FUCHSIAN HOLONOMY}
\author{LORENZO RUFFONI}
\address{Dipartimento di Matematica - Universit\`a di Bologna, Piazza di Porta San Donato 5, 40126, Bologna, Italy}
\curraddr{Department of Mathematics - Florida State University, 1017 Academic Way, Tallahassee, FL 32306-4510, USA}
\email{lorenzo.ruffoni2@gmail.com}

\keywords{Complex projective structures, quasi-Fuchsian holonomy, grafting, bubbling, movements of branch points.}
\subjclass[2010]{57M50, 20H10, 14H15}
\date{\today}

\begin{abstract}
For a given quasi-Fuchsian representation $\rho:\pi_1(S)\to$ PSL$_2\mathbb{C}$ of the fundamental group of a closed surface $S$ of genus $g\geq 2$, we prove that a generic branched complex projective structure on $S$ with holonomy $\rho$ and two branch points can be obtained from some unbranched structure on $S$ 
with the same holonomy by bubbling, i.e. a suitable connected sum with a copy of $\cp$.
\end{abstract}

\maketitle
\tableofcontents

\section{Introduction}
A complex projective structure on a surface $S$ is a geometric structure locally modelled on the geometry of the Riemann sphere $\cp$ with its group of holomorphic automorphisms $\pslc$. Since $\hyp^2, \mathbb{E}^2$ and $\mathbb{S}^2$ admit models in $\cp$, these structures generalise the classical setting of constant curvature geometries; in particular, structures with (quasi-)Fuchsian holonomy play a central role in the theory of (simultaneous) uniformization of Riemann surfaces of genus $g\geq 2$ (see \cite{P},\cite{BE}).\par
If $\rho:\pi_1(S) \to \pslr$ is a \qf representation, then the quotient of the domain of discontinuity 
of $\rho$ by the image of $\rho$ is endowed with a natural complex projective structure 
$\sigma_\rho$ with holonomy $\rho$, namely a hyperbolic structure. A natural problem is to try to 
obtain every other projective structure with holonomy $\rho$ from this hyperbolic structure via 
some elementary geometric surgeries on it. The main result in this direction is due to Goldman, who proved in 
\cite{GO} that any complex projective structure with \qf holonomy is obtained from the hyperbolic 
structure $\sigma_\rho$ via $2\pi$-grafting, i.e. by replacing some disjoint simple closed geodesic with 
Hopf annuli. As Goldman observed, a similar statement actually holds for the case of quasi-Fuchsian representations, which can be easily reduced to the Fuchsian case by a quasi-conformal map. More recently Baba has obtained analogous results for the case of Schottky representations (\cite{BA12}) and then for the generic case of purely loxodromic representations (\cite{BA15},\cite{BA17}). \par
We are interested in the analogous problem for branched complex projective structures: these are  a 
generalisation (introduced in \cite{MA}) of complex projective structures in which we allow cone 
points of angle $2\pi k$ for $ k \in \nat$. An easy way to introduce branch points is to perform a 
bubbling, i.e. to replace a simple arc with a full copy of the model space $\cp$ (see Section
\ref{s_bubbling} below for the precise definition of this surgery). The following 
question was posed by Gallo-Kapovich-Marden as Problem 12.1.2 in \cite{GKM}:\par
\vspace{.3cm}
\noindent \textbf{Question}: {\it Given two branched complex projective structures with the same holonomy, is it 
possible to pass from one to the other using the operations of grafting, degrafting, bubbling and 
debubbling?}\par
\vspace{.3cm}
Calsamiglia-Deroin-Francaviglia   provided in \cite{CDF} a positive answer in the case of 
quasi-Fuchsian holonomy, if an additional surgery is allowed, which is known as movement of branch 
points, and is a form of Schiffer variation around branch points. The main result of this paper is 
the following (see Theorem \ref{mainbubblingthm} below for the precise statement):
\begin{theorem}
Let $\rho:\pi_1(S)\to \pslc$ be quasi-Fuchsian. Then the space of branched complex projective 
structures obtained by bubbling unbranched structures with holonomy $\rho$ is 
connected, open, dense and full-measure in the moduli space of structures with two simple branch points and the 
same holonomy.
\end{theorem}
Combined with \cite[Theorem C]{GO} by Goldman this result implies a positive 
answer to the above question for a generic pair of branched complex projective structures 
with the same quasi-Fuchsian holonomy and at most two branch points. Taking into account also \cite[Theorem 
5.1]{CDF} and \cite[Theorem 1.1]{CDF2} by Calsamiglia-Deroin-Francaviglia, we can see that indeed 
the only surgeries which are generically needed to move around this moduli space are bubbling and 
debubbling. 
In the paper \cite{R} we also consider the problem of bounding the number of 
operations needed in a sequence of surgeries from one structure to another. As in Goldman's approach, for all these results the case of a general quasi-Fuchsian representation is easily reduced to that of a Fuchsian one by a quasi-conformal map.\par
The structure of the paper is the following:  Section \ref{s_bps} contains the basic definitions 
and lemmas, together with an example of a pair of non-isomorphic structures obtained by 
bubbling the same unbranched structure along isotopic arcs (see Example \ref{ex_nonisobub}); this 
phenomenon shows how sensitive these structures are to deformations. In Section 
\ref{s_geometricdecomposition} we review the geometric properties of structures with \qf 
holonomy, in the spirit of \cite{GO} and \cite{CDF}, and develop a combinatorial 
analysis of a natural decomposition of such structures into hyperbolic pieces, providing an 
explicit classification of pieces occurring for structures with at most two simple branch points; this 
already allows to prove that many structures are obtained via bubbling, and 
Section \ref{s_BMconfig} is concerned with the problem of deforming these structures without 
breaking their bubbles. Finally Section \ref{s_bubbles} contains the proof of the main theorem; the 
strategy consists in two steps: first we use the analysis in Section \ref{s_geometricdecomposition} to 
define a decomposition of the moduli space into pieces and to sort out those in which it is easy to 
find a bubbling, then we apply the results obtained in Section \ref{s_BMconfig} to move bubblings from 
these pieces to the other ones.\par
\vspace{.5cm}
\textbf{Acknowledgements}: I would like to thank Stefano Francaviglia for drawing my attention to the study of projective structures and for his constant and valuable support throughout this work. I am also very grateful to Bertrand Deroin for his interest in this project and for many useful conversations about it. Finally I thank the referees for their suggestions.

\section{Branched complex projective structures}\label{s_bps}
Let $S$ be a closed, connected and oriented surface of genus $g\geq 2$. We will denote by $\cp=\cmp 
\cup \{\infty\}$ the Riemann sphere and by $\pslc$ the group of its holomorphic automorphisms 
acting by M\"obius transformations
$$\pslc \times \cp \to \cp, \left(\begin{array}{cc} 
a & b\\ c & d \\ \end{array} \right),z \mapsto \dfrac{az+b}{cz+d}$$
We are interested in geometric structures locally modelled on this geometry, up to finite 
branched covers. The following definition is adapted from \cite{MA}.
\begin{definition}
A branched complex projective chart on $S$ is a pair $(U,\varphi)$ where $U\subset S$ is an open 
subset and $\varphi : U \to \varphi(U) \subseteq \cp$ is a finite degree orientation 
preserving branched covering map.  Two charts $(U,\varphi)$ and 
$(V,\psi)$ are compatible if $\exists \ g \in \pslc$ such that $\psi=g\varphi$ on $U\cap V$. A 
\textbf{branched complex projective structure} $\sigma$ on $S$ (BPS in the following) is the datum 
of a maximal atlas of branched complex projective charts.
\end{definition}
We will say that a structure is unbranched if all its charts are local diffeomorphisms. On the 
other hand $p \in S$ will be called a branch point of order $ord(p)=m\in \nat$ if a local chart at 
$p$ is a branched cover of degree $m+1$, i.e. if it looks like $z \mapsto z^{m+1}$.
Notice that a local chart $(U,\varphi)$ can always be shrunk to ensure 
that it contains at most one branch point and both $U$ and $\varphi(U)$ are 
homeomorphic to disks. In particular branch points are isolated, hence in finite number since 
$S$ is compact.
\begin{definition}\label{def_order}
The branching divisor of a BPS $\sigma$ is defined to be $div(\sigma)=\sum_{p\in S}ord(p)p$ and the 
branching order $ord(\sigma)$ of $\sigma$ is defined to be the degree of its branching divisor. We 
can also specify precise patterns of branching by extending this notation: for a partition 
$\lambda=(\lambda_1,\dots, \lambda_n) \in \nat^n$ we say that $\sigma$ has order 
$ord(\sigma)=\lambda$ if $div(\sigma)=\sum_{i=1}^n\lambda_i p_i$. 
\end{definition}

\begin{rmk}\label{bpsasgx}
A BPS on $S$ can be considered as a generalised ($\pslc,\cp$)-structure in the sense of 
\cite{DU}), for which the developing map may have critical points, corresponding to branch points. A 
developing map for such a structure is an orientation preserving smooth map $dev: \widetilde{S} \to 
\mathbb{CP}^1$ with isolated critical points and equivariant with respect to a holonomy 
representation $\rho:\pi_1(S)\to \pslc$. As usual, for any $g\in\pslc$ the pairs $(dev,\rho)$ and 
$(g dev, g \rho g^{-1})$ define the same BPS. Notice that in our setting $\widetilde{S}$ is a disk, 
hence $dev$ can not be a global diffeomorphism, so that these structures are never complete. Even 
worse, these structures are in general not even uniformizable, in the sense that in general $dev$ 
fails to be a diffeomorphism onto a domain $\Omega \subset \cp$ and is actually wildly non-injective. This is of course clear for branched structures, but it is actually already true in 
absence of branch points.
\end{rmk}
Let us give a few motivating examples for the study of BPSs.
\begin{example}\label{ex_brcov}
Every Riemann surface $X$ admits a non-constant meromorphic function $f$, which realizes it as a 
finite branched cover of $\cp$. This endows $X$ with a BPS with trivial holonomy and 
developing map given by $f$ itself. 
\end{example}
\begin{example}\label{ex_geom}
Every surface $S$ of genus $g\geq 2$ admits a complete Riemannian metric $g$ of constant curvature 
$-1$, which realises it as a quotient of $\mathbb{H}^2$ by a group of isometries acting freely and 
properly discontinuously. Embedding $\hyp^2$ as the upper-half plane 
$\mathcal{H}^+=\{Im(z)>0\}\subset \cmp \subset \cp$ shows that 2-dimensional hyperbolic geometry 
($\pslr,\hyp^2$) is a subgeometry of 1-dimensional complex projective geometry 
($\pslc,\mathbb{CP}^1$). Therefore these hyperbolic structures provide examples of (unbranched) 
complex projective structures. More generally it follows from the work of Troyanov in \cite{TR} 
that, given $p_1,\dots,p_n \in S$ and $k_1,\dots, k_n \in \nat$, if $\chi(S)+\sum_{i=1}^n k_i<0$ 
(resp. $=0$, or $=1$) then there exists a hyperbolic (resp. Euclidean, or spherical) metric on $S 
\setminus \{p_1,\dots,p_n \}$ with a conical singularity  of angle $2\pi (k_i+1)$ at $p_i$.
These conical hyperbolic (resp. Euclidean, or spherical) structures are examples of genuinely 
branched complex projective structures.
\end{example}
In order to define the deformation space of BPSs let us introduce a natural notion of isomorphism 
for these structures. 
\begin{definition}
 Let $\sigma$ and $\tau$ be a pair of BPSs. A map $f:\sigma\to \tau$ is projective if in local 
projective charts it is given by the restriction of a global holomorphic map 
$F:\cp \to \cp$. We say it is a projective isomorphism if it is also bijective.
\end{definition}
Recalling that any global holomorphic function $F:\cp \to \cp$ is a rational function and that the 
invertible elements in $\cmp(z)$ are exactly the fractional linear transformations 
$\frac{az+b}{cz+d}$ given by the action of $\pslc$, on obtains that a projective isomorphism is a 
diffeomorphism locally given by the restriction of some $g\in\pslc$.
\begin{definition}
A marked branched complex projective structure on $S$ is a pair $(\sigma,f)$ where $\sigma$ is a 
surface endowed with a BPS and $f:S\to \sigma$ is an orientation preserving diffeomorphism. Two 
marked BPSs $(\sigma,f)$ and $(\tau,g)$ are declared to be equivalent if $gf^{-1}:\sigma\to \tau$ is 
isotopic to a projective isomorphism $h:\sigma\to \tau$. We denote by $\mathcal{BP}(S)$ the set of 
marked branched complex projective structures on $S$ up to this equivalence relation.
\end{definition}
Thinking of BPSs in terms of equivalence classes of development-holonomy pairs $[(dev,rho)]$ as in Remark 
\ref{bpsasgx} allows us to put a natural topology (namely the compact-open topology) on this set, 
and to define a natural projection to the character variety $\chi(S)=Hom(\pi_1(S),\pslc)//\pslc$ by 
sending a BPS to its holonomy
$$hol:\mathcal{BP}(S)\to \chi(S),[\sigma]=[(dev,\rho)]\mapsto [\rho]$$
We are interested in the study of structures with a fixed holonomy, therefore we introduce the 
following subspaces of the fibres of the holonomy map.
\begin{definition}
Let $\rho\in \chi(S)$, $k \in \mathbb{N}$ and let $\lambda$ be a partition of $k$. We define
$$\mathcal{M}_{k,\rho}=\{ \sigma \in \mathcal{BP}(S) \ | \ ord(\sigma)=k, hol(\sigma)=\rho\}$$
$$\mathcal{M}_{\lambda,\rho}=\{ \sigma \in \mathcal{BP}(S) \ | \ ord(\sigma)=\lambda, 
hol(\sigma)=\rho\}$$
where the order of a structure is the one defined in Definition \ref{def_order}.
We call the \textbf{principal stratum} of $\mathcal{M}_{k,\rho}$ the subspace given by the 
partition $\lambda = (1,\dots,1)$, i.e. the one in which all branch points are simple.
\end{definition}
Recall that a representation is said to be elementary if it has a finite orbit in the standard action on $\mathbb H^3\cup \cp$, and non-elementary otherwise.
In the Appendix of \cite{CDF} Calsamiglia-Deroin-Francaviglia obtained several results about the topology of the holonomy fibre $\mathcal{M}_{k,\rho}$ (see \cite[Theorem A.2, Lemma A.7, Lemma A.13]{CDF}), which we can summarize as follows.
\begin{theorem}\label{t_complex_manifold}
Given a non-elementary representation $\rho:\pi_1(S)\to \pslc$ there exists a canonical smooth $k$-dimensional complex manifold structure on $\mathcal{M}_{k,\rho}$. Moreover the subspace determined by a partition $\lambda$ of length $n$ is a complex submanifold of dimension $n$.
\end{theorem}
In particular the principal stratum is an open dense complex 
submanifold of $\mathcal{M}_{k,\rho}$. These complex structures are locally modelled on products of 
Hurwitz spaces, i.e. spaces of deformations of finite branched cover of disks, and local 
coordinates admit a nice geometric description (see Remark \ref{movingiscoord} below for more details about description of local neighbourhoods for this manifold structure.).
\begin{rmk}\label{uniquedev}
In the following, when working with a BPS $\sigma$, we will find it convenient to fix a 
representative representation $\rho$ of the holonomy $hol(\sigma)$, i.e. to choose a representation 
in its conjugacy class. As soon as the holonomy is non-elementary, there will be a unique developing 
map equivariant with respect to the chosen representation. Indeed if $dev_1$ and $dev_2$ are 
developing maps for $\sigma$ equivariant with respect to $\rho:\pi_1(S)\to \pslc$, then $\exists \ 
g \in \pslc$ such that $dev_2=g dev_1$ and for any $\gamma\in \pi_1(S)$  we have
$$\rho(\gamma)g dev_1 = \rho(\gamma) dev_2 = dev_2 \gamma =g dev_1 \gamma =g \rho(\gamma) dev_1$$
so that $(\rho(\gamma)g)^{-1}g \rho(\gamma)$ is an element of $\pslc$ fixing every point of 
$dev_1(\widetilde{S})$. Since a developing map has isolated critical points, there is some point 
of $\widetilde{S}$ at which it is a local diffeomorphism, hence its image has non-empty interior. 
But a  M\"obius transformation fixing more than three points is the identity of $\cp$, hence  
$(\rho(\gamma)g)^{-1}g \rho(\gamma)=id$. This means that $g$ is in the centralizer of the image of 
$\rho$, which is trivial since the holonomy is assumed to be non-elementary; as a result
$g=id$ and the two developing maps coincide.
\end{rmk}
We conclude this preliminary section by introducing three elementary geometric surgeries which one 
can perform on a given BPS to obtain a new BPS with the same holonomy.

\subsection{Grafting}\label{s_grafting}
The first surgery consists in replacing a simple closed curve with an annulus endowed with a 
projective structure determined by the structure we begin with. It was first introduced by Maskit in 
\cite{MAS2} to produce examples of projective structures with surjective developing map; here we 
review it mainly to fix terminology and notation. Let us pick  $\sigma \in \mathcal{BP}(S)$ and let 
$(dev,\rho)$ be a development-holonomy pair defining it.
\begin{definition}
 Let $\gamma \subset S$ be a simple closed curve on $S$. We say that $\gamma$ is graftable with 
respect to $\sigma$ if $\rho(\gamma)$ is loxodromic (i.e. not elliptic nor parabolic) and $\gamma$ 
is injectively developed, i.e. the 
restriction of $dev$ to any of its lifts $\widetilde{\gamma}\subset \widetilde{S}$ is 
injective.
\end{definition}
Since $dev$ is $\rho$-equivariant, if $\gamma$ is graftable then a developed image of it is an 
embedded arc in $\cp$ joining the two fixed points of $\rho(\gamma)$. Moreover $\rho(\gamma)$ acts 
freely and properly discontinuously on $\cp \setminus \overline{dev(\widetilde{\gamma})}$ and the 
quotient is an annulus endowed with a complete unbranched complex projective structure. 
\begin{definition}
 Let $\gamma\subset S$ be a graftable curve with respect to $\sigma$. For any lift 
$\widetilde{\gamma}$ of $\gamma$ we cut $\widetilde{S}$ along it and a copy of $\cp$ along 
$\overline{dev(\widetilde{\gamma})}$, and glue them together equivariantly via the developing map. 
This gives us a simply connected surface $\widetilde{S}'$ to which the action 
$\pi_1(S)\curvearrowright \widetilde{S}$ and the map $dev:\widetilde{S}\to \cp$ naturally extend, 
so that the quotient gives rise to a new structure $\sigma ' \in \mathcal{BP}(S)$. We call this 
structure the \textbf{grafting} of $\sigma$ along $\gamma$ and denote it by $Gr(\sigma,\gamma)$.
The surface $\sigma \setminus \gamma$ projectively embeds in $Gr(\sigma, \gamma)$ and  the 
complement is the annulus $A_\gamma=(\cp \setminus\overline{dev(\widetilde{\gamma})})/\rho(\gamma)$, 
which we call the grafting annulus associated to $\gamma$.  The inverse operation will be called a 
degrafting.
\end{definition}
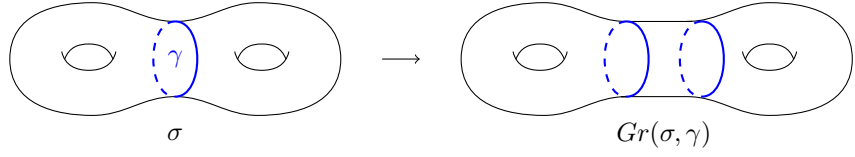
\begin{figure}[h]

\begin{center}

\begin{tikzpicture}[scale=0.5]
\node at (0,-2) {$\sigma$};
\node at (0,0) {\color{blue} $\gamma$};
\draw[xscale=0.8] (-5.5,0) to[out=90,in=180] (-2.75,1.5) to[out=0,in=180] (0,1) to[out=0,in=180] (2.75,1.5) to[out=0,in=90] (5.5,0) ;
\draw[rotate=180,xscale=0.8] (-5.5,0) to[out=90,in=180] (-2.75,1.5) to[out=0,in=180] (0,1) to[out=0,in=180] (2.75,1.5) to[out=0,in=90] (5.5,0) ;

\draw[xshift=-0.3cm,yshift=-0.1cm,xscale=0.8]  (-3.3,0.1) to[out=65,in=180] (-2.5,0.5) to[out=0,in=115] (-1.7,0.1);
\draw[xshift=-0.3cm,yshift=-0.1cm,xscale=0.8]  (-3.4,0.3) to[out=-75,in=180] (-2.5,-0.2) to[out=0,in=-105] (-1.6,0.3);
\draw[xshift=4.3cm,yshift=-0.1cm,xscale=0.8]  (-3.3,0.1) to[out=65,in=180] (-2.5,0.5) to[out=0,in=115] (-1.7,0.1);
\draw[xshift=4.3cm,yshift=-0.1cm,xscale=0.8]  (-3.4,0.3) to[out=-75,in=180] (-2.5,-0.2) to[out=0,in=-105] (-1.6,0.3);

\draw[blue,dashed,thick] (0,1) to[out=180,in=180] (0,-1);
\draw[blue,thick] (0,1) to[out=0,in=0] (0,-1);


\draw[->] (5.5,0) -- (6.5,0);

\begin{scope}[xshift=12cm]
 \node at (1,-2) {$Gr(\sigma,\gamma)$};
\draw[xscale=0.8] (-5.5,0) to[out=90,in=180] (-2.75,1.5) to[out=0,in=180] (0,1) to[out=0,in=0] (2,1) 
to[out=0,in=180] (4.75,1.5) to[out=0,in=90] (7.5,0) ;
\draw[rotate=180,xscale=0.8,xshift=-2cm] (-5.5,0) to[out=90,in=180] (-2.75,1.5) to[out=0,in=180] 
(0,1) to[out=0,in=0] (2,1) to[out=0,in=180] (4.75,1.5) to[out=0,in=90] (7.5,0) ;

\draw[xshift=-0.3cm,yshift=-0.1cm,xscale=0.8]  (-3.3,0.1) to[out=65,in=180] (-2.5,0.5) 
to[out=0,in=115] (-1.7,0.1);
\draw[xshift=-0.3cm,yshift=-0.1cm,xscale=0.8]  (-3.4,0.3) to[out=-75,in=180] (-2.5,-0.2) 
to[out=0,in=-105] (-1.6,0.3);
\draw[xshift=5.8cm,yshift=-0.1cm,xscale=0.8]  (-3.3,0.1) to[out=65,in=180] (-2.5,0.5) 
to[out=0,in=115] (-1.7,0.1);
\draw[xshift=5.8cm,yshift=-0.1cm,xscale=0.8]  (-3.4,0.3) to[out=-75,in=180] (-2.5,-0.2) 
to[out=0,in=-105] (-1.6,0.3);

\draw[blue,dashed,thick] (0,1) to[out=180,in=180] (0,-1);
\draw[blue,thick] (0,1) to[out=0,in=0] (0,-1);
\draw[blue,dashed,thick] (2,1) to[out=180,in=180] (2,-1);
\draw[blue,thick] (2,1) to[out=0,in=0] (2,-1);
\end{scope}

\end{tikzpicture}

\end{center}
\caption{Grafting a surface}
\end{figure}

The easiest example of this construction consists in grafting a simple geodesic on a hyperbolic 
surface; for such a structure every simple essential curve $\gamma$ is graftable, since 
the holonomy is purely hyperbolic and the developing map is globally injective. 
The grafting surgery preserves the holonomy and does not involve any modification of the branching 
divisor, so that if $\sigma \in \mathcal{M}_{\lambda,\rho}$ then $Gr(\sigma,\gamma) 
\in \mathcal{M}_{\lambda,\rho}$ too. Notice that for any structure $\sigma$ and any 
graftable curve $\gamma$ on it the structure $Gr(\sigma, \gamma)$ has surjective but non-injective 
developing map.

\subsection{Bubbling}\label{s_bubbling}
The second surgery consists in replacing a simple arc with a disk endowed with a projective 
structure determined by the structure we begin with, hence it can be thought as a ``finite version''
of grafting. It was first considered by Gallo-Kapovich-Marden  in \cite{GKM} as a tool to introduce 
new branch points on a projective structures. As before, let us choose  $\sigma \in 
\mathcal{BP}(S)$ and let $(dev,\rho)$ be a development-holonomy pair defining it.
\begin{definition}
 Let $\beta \subset S$ be a simple arc on $S$. We say that $\beta$ is bubbleable with 
respect to $\sigma$ if it is injectively developed, i.e. the restriction of $dev$ to any of its 
lifts $\widetilde{\beta}\subset \widetilde{S}$ is injective.
\end{definition}
The surgery is then defined as follows.
\begin{definition}
 Let $\beta\subset S$ be a bubbleable arc with respect to $\sigma$. For any lift 
$\widetilde{\beta}$ of $\beta$ we cut $\widetilde{S}$ along it and a copy of $\cp$ along 
$dev(\widetilde{\beta})$, and glue them together equivariantly via the developing map. 
Once again, this gives us a simply connected surface $\widetilde{S}'$ to which the action 
$\pi_1(S)\curvearrowright \widetilde{S}$ and the map $dev:\widetilde{S}\to \cp$ naturally extend, 
so that the quotient gives rise to a new structure $\sigma' \in \mathcal{BP}(S)$. We call this 
structure the \textbf{bubbling} of $\sigma$ along $\beta$ and denote it by $Bub(\sigma,\beta)$.
The surface $\sigma \setminus \beta$ projectively embeds in $Bub(\sigma,\beta)$ and the 
complement is the disk $B=\cp  \setminus dev(\widetilde{\beta})$, which we call the bubble 
associated to $\beta$. 
\end{definition}
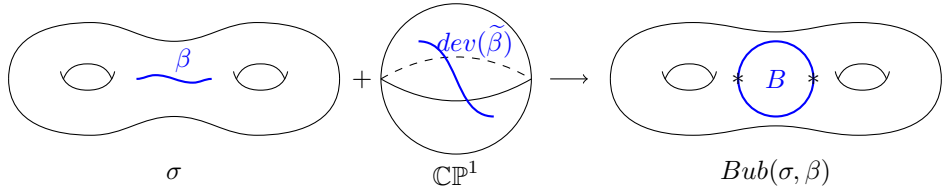
\begin{figure}[h]

\begin{center}
\begin{tikzpicture}[scale=0.5]
\node at (0,-2.5) {$\sigma$};
\draw[xscale=0.8] (-5.5,0) to[out=90,in=180] (-2.75,1.5) to[out=0,in=180] (0,1) to[out=0,in=180] (2.75,1.5) to[out=0,in=90] (5.5,0) ;
\draw[rotate=180,xscale=0.8] (-5.5,0) to[out=90,in=180] (-2.75,1.5) to[out=0,in=180] (0,1) to[out=0,in=180] (2.75,1.5) to[out=0,in=90] (5.5,0) ;

\draw[xshift=-0.3cm,yshift=-0.1cm,xscale=0.8]  (-3.3,0.1) to[out=65,in=180] (-2.5,0.5) to[out=0,in=115] (-1.7,0.1);
\draw[xshift=-0.3cm,yshift=-0.1cm,xscale=0.8]  (-3.4,0.3) to[out=-75,in=180] (-2.5,-0.2) to[out=0,in=-105] (-1.6,0.3);
\draw[xshift=4.3cm,yshift=-0.1cm,xscale=0.8]  (-3.3,0.1) to[out=65,in=180] (-2.5,0.5) to[out=0,in=115] (-1.7,0.1);
\draw[xshift=4.3cm,yshift=-0.1cm,xscale=0.8]  (-3.4,0.3) to[out=-75,in=180] (-2.5,-0.2) to[out=0,in=-105] (-1.6,0.3);

\node at (0.25,0.5) {\color{blue} $\beta$};
\draw[blue,thick] (-1,0) to[out=0,in=180] (-0.5,0.1) to[out=0,in=180] (0.5,-0.1) to[out=0,in=180] (1,0);

\node at (5,0) {$+$};

\begin{scope}[xshift=-2.5cm]
\node at (10,-2.5) {$\mathbb{CP}^1$};
\node at (10.5,1) {\color{blue} $dev(\widetilde{\beta})$};
\draw (8,0) to[out=90,in=180] (10,2) to[out=0,in=90] (12,0) to[out=-90,in=0] (10,-2)  to[out=180,in=-90] (8,0);
\draw (8,0) to[out=-30,in=210] (12,0);
\draw[dashed] (8,0) to[out=30,in=150] (12,0);
\draw[blue,thick] (9,1) to[out=0,in=180] (11,-1);
\draw[->] (12.5,0) -- (13.5,0);
\end{scope}

\begin{scope}[xshift=16cm]
\node at (0,-2.5) {$Bub(\sigma,\beta)$};
\draw[xscale=0.8] (-5.5,0) to[out=90,in=180] (-2.75,1.5) to[out=0,in=180] (0,1.25) to[out=0,in=180] (2.75,1.5) to[out=0,in=90] (5.5,0) ;
\draw[rotate=180,xscale=0.8] (-5.5,0) to[out=90,in=180] (-2.75,1.5) to[out=0,in=180] (0,1.25) to[out=0,in=180] (2.75,1.5) to[out=0,in=90] (5.5,0) ;

\draw[xshift=-0.3cm,yshift=-0.1cm,xscale=0.8]  (-3.3,0.1) to[out=65,in=180] (-2.5,0.5) to[out=0,in=115] (-1.7,0.1);
\draw[xshift=-0.3cm,yshift=-0.1cm,xscale=0.8]  (-3.4,0.3) to[out=-75,in=180] (-2.5,-0.2) to[out=0,in=-105] (-1.6,0.3);
\draw[xshift=4.3cm,yshift=-0.1cm,xscale=0.8]  (-3.3,0.1) to[out=65,in=180] (-2.5,0.5) to[out=0,in=115] (-1.7,0.1);
\draw[xshift=4.3cm,yshift=-0.1cm,xscale=0.8]  (-3.4,0.3) to[out=-75,in=180] (-2.5,-0.2) to[out=0,in=-105] (-1.6,0.3);

\draw[blue,thick] (-1,0) to[out=90,in=180] (0,1) to[out=0,in=90] (1,0);
\draw[blue,thick] (-1,0) to[out=-90,in=180] (0,-1) to[out=0,in=-90] (1,0);
\node at (-1,0) {$*$};
\node at (1,0) {$*$};
\node at (0,0) {\color{blue} $B$};

\end{scope}

\end{tikzpicture}

\end{center}
\caption{Bubbling a surface}
\end{figure}
The easiest example is obtained by bubbling a hyperbolic surface along an embedded geodesic arc.
The bubbling surgery preserves the holonomy and introduces a pair of simple branch points 
corresponding to the endpoints of the bubbling arc. Therefore if $\sigma \in 
\mathcal{M}_{\lambda,\rho}$ then $Bub(\sigma,\beta) \in \mathcal{M}_{\lambda+(1,1),\rho}$, where if 
$\lambda$ is a partition of $k$, $\lambda+(1,1)$ is the partition of $k+2$ obtained appending 
$(1,1)$ to it.\par
Once  a bubbling is performed, we see a subsurface of $S$ homeomorphic to a disk and  
isomorphic to $\mathbb{CP}^1$ cut along a simple arc, the isomorphism being given by any 
determination of the developing map itself. It is useful to be able to recognise this kind of 
subsurface, since there is an obvious way to remove it and lower the branching order by $2$; such 
operation is called debubbling and is the inverse of bubbling. 
\begin{definition}
A \textbf{bubble} on $\sigma \in \mathcal{BP}(S)$ is an embedded closed disk $B \subset S$ 
 whose boundary decomposes as $\partial B=\beta' \cup \{x,y\} \cup \beta''$ where $\{x,y\}$ are 
simple branch points of $\sigma$ and $\beta',\beta''$ are embedded injectively developed arcs which 
overlap once developed; more precisely there exist a determination of the developing map on $B$ 
which injectively maps $\beta',\beta''$ to the same simple arc $\widehat{\beta}\subset \cp$ and 
restricts to a diffeomorphism $dev: int(B)\to \cp \setminus \widehat{\beta} $.
\end{definition}
\begin{figure}[h]
\input{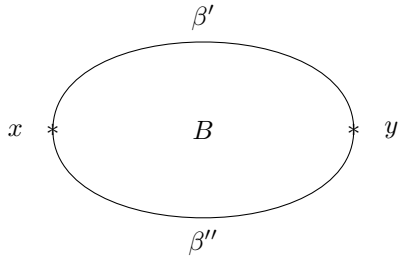}
\caption{A bubble}
\end{figure}
Notice that a BPS obtained by bubblings some unbranched structures has by definition an even 
number of branch points and surjective non-injective developing map. As a consequence branched 
hyperbolic structures with an even number of branch points do not arise as bubblings, as their 
developing maps take value only in the upper-half plane. By the work of \cite{TR} these structures 
exist on every surface of genus $g\geq 3$; this example was already mentioned in \cite{CDF}.

\subsection{Movements of branch points}\label{s_movementsofbranchpoints}
The last surgery we will use takes place locally around a branch point and consists in a 
deformation of the local branched projective chart, which can be thought as an analogue in our 
setting of the Schiffer variations in the theory of Riemann surfaces (see \cite{N}). They were 
introduced by Tan in \cite{TA} for simple branch points (and then generalised in \cite{CDF} for 
branch points of higher order) as a tool to perform local deformations of a BPS inside the moduli 
space $\mathcal{M}_{k,\rho}$. Since we will need this surgery only for simple branch points, we 
restrict here to that case 
and avoid the  technicalities required by a more general treatment.
\begin{definition}
Let $\sigma \in \mathcal{M}_{k,\rho}$ and let $p\in \sigma$ be a simple branch point. An 
\textbf{embedded twin pair} at $p$ is a pair of embedded arcs $\mu=\{\mu_1,\mu_2\}$ which 
meet exactly at $p$, are injectively developed and overlap once developed; more precisely there 
exist a determination of the developing map around $\mu_1 \cup \mu_2$ which injectively maps 
$\mu_1,\mu_2$ to the same simple arc $\widehat{\mu}\subset \cp$. 
\end{definition}
Given such a pair of arcs we can perform the following cut-and-paste surgery.
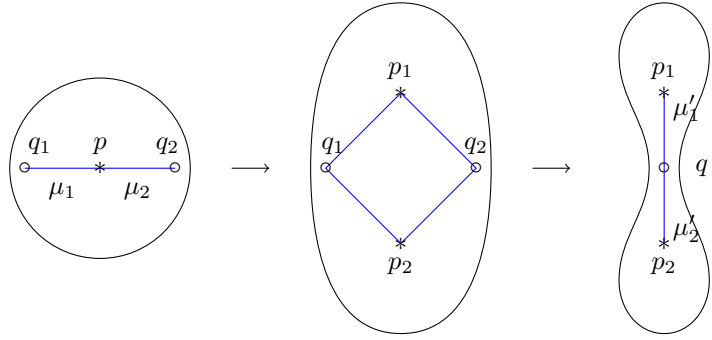
\begin{figure}[h]

\begin{center}
\begin{tikzpicture}

\draw[xshift=-4cm] (-1.2,0) to[out=90,in=180] (0,1.2) to[out=0,in=90] (1.2,0) to[out=-90,in=0] (0,-1.2) to[out=180,in=-90] (-1.2,0);
\draw[xshift=-4cm,blue] (-1,0) to (1,0);
\node[xshift=-4cm] at (-0.5,-0.3) {$\mu_1$};
\node[xshift=-4cm] at (0.5,-0.3) {$\mu_2$};
\node[xshift=-4cm] at (0,0) {$*$};
\node[xshift=-4cm] at (0,0.3) {$p$};
\node[xshift=-4cm] at (1,0) {$\circ$};
\node[xshift=-4cm] at (0.9,0.3) {$q_2$};
\node[xshift=-4cm] at (-1,0) {$\circ$};
\node[xshift=-4cm] at (-0.8,0.3) {$q_1$};

\draw[->] (-2.25,0) to (-1.75,0);

\draw (-1.2,0) to[out=90,in=180] (0,2.2) to[out=0,in=90] (1.2,0) to[out=-90,in=0] (0,-2.2) to[out=180,in=-90] (-1.2,0);
\draw[blue] (-1,0) to (0,1) to (1,0) to (0,-1) to (-1,0);
\node at (0,1) {$*$};\node at (0,1.3) {$p_1$};
\node at (0,-1) {$*$};\node at (0,-1.3) {$p_2$};
\node at (1,0) {$\circ$};\node at (1,0.3) {$q_2$};
\node at (-1,0) {$\circ$};\node at (-0.9,0.3) {$q_1$};

\draw[->] (1.75,0) to (2.25,0);

\draw[xshift=3.5cm] (-0.2,0) to[out=90,in=-90] (-0.6,1.5) to[out=90,in=180] (0,2.2) to[out=0,in=90] (0.6,1.5) to[out=-90,in=90] (0.2,0) to[out=-90,in=90] (0.6,-1.5) to[out=-90,in=0] (0,-2.2) to[out=-180,in=-90] (-0.6,-1.5) to[out=90,in=-90] (-0.2,0);
\draw[xshift=3.5cm,blue] (0,1) to (0,-1);\node[xshift=3.5cm] at (0.3,0.8) 
{$\mu'_1$};\node[xshift=3.5cm] at (0.3,-0.8) {$\mu'_2$};
\node[xshift=3.5cm] at (0,1) {$*$};\node[xshift=3.5cm] at (0,1.3) {$p_1$};
\node[xshift=3.5cm] at (0,-1) {$*$};\node[xshift=3.5cm] at (0,-1.3) {$p_2$};
\node[xshift=3.5cm] at (0,0) {$\circ$};\node[xshift=3.5cm] at (0.5,0) {$q$};

\end{tikzpicture}
\end{center}
\caption{A movement of branch point}\label{pic_schiffer}
\end{figure}
\begin{definition}
Let $\sigma \in \mathcal{M}_{k,\rho}$, let $p\in \sigma$ be a simple branch point and 
$\mu=\{\mu_1,\mu_2\}$ an embedded twin pair at $p$. The BPS $\sigma'$ obtained by cutting $S$ along 
$\mu_1\cup \mu_2$ and regluing the resulting boundary with the obvious identification (as shown in 
Figure \ref{pic_schiffer}) is said to be obtained by a \textbf{movement of branch point} at $p$ 
along $\mu$ and is denoted by $Move(\sigma,\mu)$.
\end{definition}
This surgery preserves the holonomy and does not change the structure of the branching divisor. 
Therefore if $\sigma \in \mathcal{M}_{\lambda,\rho}$ then $Move(\sigma,\mu) \in 
\mathcal{M}_{\lambda,\rho}$. Notice that the image of the developing map is not changed by this 
operation; moreover once a movement is performed, we have an induced embedded twin pair on the new 
structure, and moving points along it of course brings us back to $\sigma$.
\begin{rmk}
The movement of branch points along an embedded twin pair $\mu$ is a deformation which comes in a 
1-parameter family. In the above notations, if $\widehat{\mu}:[0,1]\to \cp$ is a parametrization of 
the developed image of $\mu$, then for $t\in[0,1]$ we can consider the structure 
$\sigma_t=Move(\sigma,\mu^t)$, where $\mu^t$ is the embedded twin pair contained in $\mu$ and 
developing to the subarc 
$\widehat{\mu}([0,t])$. Following this deformation as $t$ varies, we see the developed image of the 
branch point sliding along the arc $\widehat{\mu}$, which motivates the name of this surgery.
\end{rmk}
\begin{rmk}\label{movingiscoord}
As anticipated above (see Theorem \ref{t_complex_manifold}), as soon as the holonomy $\rho$ is non-elementary, the moduli space 
$\mathcal{M}_{\lambda,\rho}$ carries a natural structure of complex manifold of dimension equal to 
the length of the partition $\lambda$. It is proved in the Appendix of \cite{CDF} that the local 
neighbourhoods of a BPS $\sigma \in \mathcal{M}_{k,\rho}$ for this topology are obtained by local 
deformations at the branch points. Just by counting dimensions, we see that for simple branch points (i.e. structures in the principal 
stratum) these are just the movements of branch points described above; for higher order 
branch points one needs to introduce a slight generalisation of them, but we will not need this.
\end{rmk}

\subsection{Injectively developed isotopies}
We have so far introduced some surgeries which can be performed on a BPS, which depend on the 
choice of a simple arc which is injectively mapped to $\cp$ by the developing map. It is natural to 
ask how much this choice is relevant as far as the isomorphism class of the resulting structure is 
concerned; an answer to this will be needed in the forthcoming sections. The following turns out 
to be the a useful notion to consider.
\begin{definition}
 Let $\sigma$ be a BPS on $S$ and $\eta:[0,1]\to S$ an embedded arc with embedded developed image. 
An isotopy $H:[0,1]\times [0,1]\to S$ of $\eta$ is said to be injectively developed if 
$\eta_s=H(s,.)$ is an embedded arc with embedded developed image for all $s \in [0,1]$.
\end{definition}
Let us begin with the following lemma, which says that being injectively developed is a stable 
condition.
\begin{lemma}\label{injdevnbdofinjdevpath}
 Let $\sigma$ be a BPS on $S$. Let $\gamma:[0,1]\to S$ be an embedded arc having embedded developed 
image and not going through branch points (except possibly at its endpoints). Then there exists an 
injectively developed subset $U\subseteq S$ such that $\gamma \subset U$ and $\gamma(]0,1[) \subset 
int(U)$.
\end{lemma}
\begin{proof}
Let  $\widetilde{\gamma}$ be a lift of the arc to the universal cover. 
Assume first that $\gamma$ does not go through any branch point at all. Then we can prove that it 
actually has an injectively developed neighbourhood: if this were not the case, there would be a 
sequence of nested open neighbourhoods $U_{n+1}\subsetneq U_n$ of $\gamma$ such that $\forall \ n 
\in \nat$ we could find a pair of distinct point $x_n,y_n \in U_n$  with the same developed 
image. By compactness of $S$, these sequences subconverge to a pair of points $x,y\in \gamma$ 
with the same developed image. Since the path is injectively developed, we get $x=y$. But since the 
path does not go through branch points, the developing map is locally injective at any of its 
points, so that the existence of the points $x_n,y_n$ arbitrarily close to $x=y\in \gamma$ is 
absurd.\\
If one endpoint, say  $\gamma(0)$, of $\gamma$ is a branch point of order $k$, then 
clearly every set containing it in its interior is not injectively developed. Nevertheless a 
sufficiently small neighbourhood $\Omega$ of $\gamma(0)$ decomposes as a 
disjoint union of injectively developed sectors $A_1,\dots,A_{k+1}$; an initial segment of 
$\gamma$ belongs to one of them, say  $A_1$; so we can simply pick a sequence of 
nested 
sets $V_{n+1}\subsetneq V_n$ such that for every $n\in \nat$ we have 
that $\gamma(0)\in V_n$, $V_n \cap \Omega \subsetneq A_1$ and $V_n$ contains 
$\gamma(]0,1[)$ in its interior, and apply the previous argument to obtain a 
pair of sequences $x_n\neq y_n \in V_n$ converging to $x=y\in \gamma$. The non trivial case to 
discuss is the case in which the limit is a branch point, i.e. $x=y=\gamma(0)$; by construction of 
$V_n$, for $n$ large enough the points $x_n,y_n$ must lie inside $A_1$, which is injectively 
developed, hence we reach a contradiction exactly as before.
\end{proof}

In particular this implies that it is always possible to perform small deformations of an 
injectively 
developed arc through an injectively developed isotopy relative to endpoints. This applies both to 
bubbleable arcs and to arcs appearing in an embedded twin pair. Injectively developed 
isotopies of bubbleable arcs and embedded twin pairs are relevant in our discussion since they do 
not change the isomorphism class of the structure obtained by performing a bubbling or a movement 
of branch points, as established by the following statements. The next one is simply a 
reformulation of \cite[Lemma 2.8]{CDF}.

\begin{lemma}\label{bubblingonisotopicarcsareiso} 
Let $\sigma$ be a BPS and let $\beta,\beta' \subset \sigma$ be bubbleable arcs with the same 
endpoints. 
If there exists an injectively developed isotopy relative to endpoints from $\beta$ to $\beta'$, 
then $Bub(\sigma,\beta)=Bub(\sigma,\beta')$. 
\end{lemma}

The following is the statement, analogous to Lemma \ref{bubblingonisotopicarcsareiso}, for a movement of 
branch points along different embedded twin pairs.
\begin{lemma}\label{moveisotopic}
 Let $\sigma $ be a BPS and let $p$ be a simple branch point. Let 
$\mu=\{\mu_1,\mu_2\}$ and  $\nu=\{\nu_1,\nu_2\}$ be embedded twin pairs based at $p$ with the 
same endpoints, and let $q_i$ be the common endpoint of $\mu_i$ and $\nu_i$ for $i=1,2$.
Suppose that there exists an injectively developed isotopy $H:[0,1]\times [-1,1]\to S$ from $\mu$ 
to 
$\nu$ relative to $\{q_1,p,q_2\}$ and such that $\alpha^s=\{\alpha^s_1=H(s,[-1,0])$, 
$\alpha^s_2=H(s,[0,1])\}$ is an embedded twin pair for all $s \in  [0,1]$.
Then $Move(\sigma,\mu)=Move(\sigma,\nu)$.
\end{lemma}
\begin{proof}
First of all notice that each path $\alpha^s_i$ appearing in an embedded twin pair $\alpha^s$ is in 
particular an embedded arc which is injectively developed and goes through exactly one branch 
point, 
which is $p$. Therefore we can pick an injectively developed set $U^s_i$ containing 
$\alpha^s_i\setminus \{p\}$ in its interior as in Lemma \ref{injdevnbdofinjdevpath}. We can choose 
this set in such a way that $U^s=U^s_1\cup U^s_2$ is an open neighbourhood of $\alpha^s$: for 
instance we can take $U^s_1$ such that its developed image is an open neighbourhood of the 
developed image of $\alpha^s$, then pull it back via the developing map, so that $U^s$ is the 
domain of a local projective chart which simply branches at $p$ and contains the whole embedded 
twin pair $\alpha^s$.
The sets $U^s$ provide an open cover of $Im(H)$; by compactness we extract a finite subcover 
indexed by some $s_0=0,s_1,\dots,s_N=1$. Up to taking an intermediate finite subcover between 
$\{U^{s_0},\dots,U^{s_N}\}$ and $\{U^s \ | \ s\in [0,1] \}$ we can assume that the local chart 
$U^{s_i}$ contains not only $\alpha^{s_i}$ but also $\alpha^{s_{i\pm 1}}$. Then we conclude by 
observing that $\alpha^{s_0}=\mu$ and $\alpha^{s_N}=\nu$ and that the results 
in the Appendix of \cite{CDF} (see also Theorem \ref{t_complex_manifold} and Remark \ref{movingiscoord} above) imply that $Move(\sigma, \alpha^{s_i})=Move(\sigma,\alpha^{s_{i+ 
1}})$, because $\alpha^{s_i}$ and $\alpha^{s_{i+1}}$ are contained in the domain of a single local 
chart.
\end{proof}

We conclude this preliminary section by remarking that an ordinary isotopy is in general not enough 
to obtain this kind of results. In the next example we provide an explicit construction of two 
bubbleable arcs which are isotopic but not isotopic through an  injectively developed isotopy, for 
which the resulting structures are not isomorphic. Most of the technical parts in Section 
\ref{s_BMconfig} 
below are needed to avoid this kind of phenomenon, which was already observed in \cite[Remark 
3.4]{CDF2} for the case of graftings.

\begin{example}\label{ex_nonisobub}
Let $S$ be a genus 2 surface with a hyperbolic structure $\sigma_\rho$, with holonomy a Fuchsian 
representation $\rho:\pi_1(S)\to \pslr$, and let $\gamma$ 
be a separating oriented closed geodesic. Let $\eta$ be an oriented embedded geodesic arc on $S$ 
with one endpoint $x$ on $\gamma$ and orthogonally intersecting $\gamma$ only in $x$; let $y$ be 
the 
other endpoint, which we assume to be on the right of $\gamma$ (see Figure \ref{pic_nonisobub1}). 
We want to perform a grafting of 
$\sigma_\rho$ along $\gamma$ and then show how to perform two different bubbling on 
$Gr(\sigma_\rho,\gamma)$ along two different extensions of $\eta$.
On $Gr(\sigma_\rho,\gamma)$ we have two distinguished curves $\gamma^\pm$ coming from $\gamma$ and 
bounding the grafting annulus $A_\gamma$. We also have two marked points $x^\pm \in \gamma^\pm$ 
coming from the point $x$, and an arc coming from $\eta$, which we still denote by the same name, 
which starts at $x^+ \in \gamma^+$ orthogonally and moves away from the annulus.\par
There is a natural way to extend $\eta$ by analytic continuation to an embedded arc reaching the 
other point $x^- \in \gamma^-$: namely consider the extension of the developed image of $\eta$ 
(which is a small geodesic arc in the upper half-plane) to a great circle $\widehat{\eta}$ on 
$\cp$. This gives an embedded arc on $Gr(\sigma_\rho,\gamma)$ which is not injectively 
developed, hence not bubbleable.
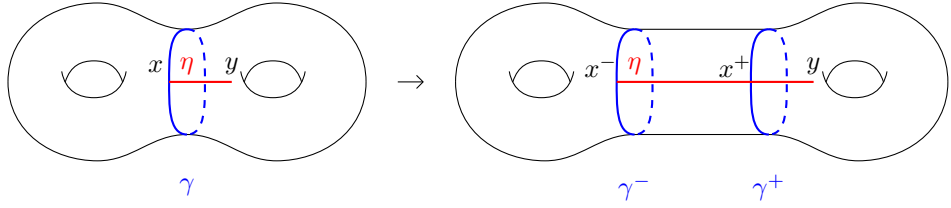
\begin{figure}[h]

\begin{center}
\begin{tikzpicture}[scale=0.7,xscale=0.85]
\draw (0,0) to[out=90,in=180] (2,1.5) to[out=0,in=180] (4,1) to[out=0,in=180] (6,1.5) 
to[out=0,in=90] (8,0) ;
\draw[yscale=-1] (0,0) to[out=90,in=180] (2,1.5) to[out=0,in=180] (4,1) to[out=0,in=180] (6,1.5) 
to[out=0,in=90] (8,0) ; 
\draw[xscale=0.8]  (1.6,0) to[out=65,in=180] (2.4,0.4) to[out=0,in=115] (3.2,0);
\draw[xscale=0.8]  (1.5,0.2) to[out=-75,in=180] (2.4,-0.3) to[out=0,in=-105] (3.3,0.2);
\draw[xshift=4cm,xscale=0.8]  (1.6,0) to[out=65,in=180] (2.4,0.4) to[out=0,in=115] (3.2,0);
\draw[xshift=4cm,xscale=0.8]  (1.5,0.2) to[out=-75,in=180] (2.4,-0.3) to[out=0,in=-105] (3.3,0.2);

\node at (4,-2) {\color{blue} $\gamma$};
\draw[blue,thick] (4,1) to[out=180,in=90] (3.6,0) to[out=-90,in=180] (4,-1);
\draw[blue,thick,dashed] (4,1) to[out=0,in=90] (4.4,0) to[out=-90,in=0] (4,-1);

\node (P) at (8.5,0) {};
 \node (R) at (9.5,0) {};
 \path[->,font=\scriptsize,>=angle 90] (P) edge node[above]{} (R);
 
\draw[xshift=10cm] (0,0) to[out=90,in=180] (2,1.5) to[out=0,in=180] (4,1) to[out=0,in=180] (7,1) 
to[out=0,in=180] (9,1.5) to[out=0,in=90] (11,0) ;
\draw[xshift=10cm,yscale=-1] (0,0) (0,0) to[out=90,in=180] (2,1.5) to[out=0,in=180] (4,1) 
to[out=0,in=180] (7,1) to[out=0,in=180] (9,1.5) to[out=0,in=90] (11,0) ;
\draw[xshift=10cm,xscale=0.8]  (1.6,0) to[out=65,in=180] (2.4,0.4) to[out=0,in=115] (3.2,0);
\draw[xshift=10cm,xscale=0.8]  (1.5,0.2) to[out=-75,in=180] (2.4,-0.3) to[out=0,in=-105] (3.3,0.2);
\draw[xshift=17cm,xscale=0.8]  (1.6,0) to[out=65,in=180] (2.4,0.4) to[out=0,in=115] (3.2,0);
\draw[xshift=17cm,xscale=0.8]  (1.5,0.2) to[out=-75,in=180] (2.4,-0.3) to[out=0,in=-105] (3.3,0.2);

\draw[xshift=10cm,blue,thick] (4,1) to[out=180,in=90] (3.6,0) to[out=-90,in=180] (4,-1);
\draw[xshift=10cm,blue,thick,dashed] (4,1) to[out=0,in=90] (4.4,0) to[out=-90,in=0] (4,-1);
\draw[xshift=13cm,blue,thick] (4,1) to[out=180,in=90] (3.6,0) to[out=-90,in=180] (4,-1);
\draw[xshift=13cm,blue,thick,dashed] (4,1) to[out=0,in=90] (4.4,0) to[out=-90,in=0] (4,-1);

\node at (14,-2) {\color{blue} $\gamma^-$};
\node at (17,-2) {\color{blue} $\gamma^+$};

\draw[red,thick] (3.6,0) to[out=0,in=180] (5,0);
\node at (3.3,0.3) { $x$};
\node at (5,0.3) { $y$};
\node at (4,0.3) {\color{red} $\eta$};
\draw[red,thick] (13.6,0) to[out=0,in=180] (18,0);
\node at (13.25,0.3) { $x^-$};
\node at (16.25,0.3) { $x^+$};
\node at (18,0.3) { $y$};
\node at (14,0.3) {\color{red} $\eta$};

\end{tikzpicture}

\end{center}
\caption{Analytic extension of $\eta$ in $Gr(\sigma_\rho,\gamma)$}\label{pic_nonisobub1}
\end{figure}
To obtain bubbleable arcs we slightly perturb this construction; in $\cp$ consider an embedded arc 
which starts at the developed image $\widehat{x}$ of $x$ and ends at the developed image 
$\widehat{y}$ of $y$, but leaves $\widehat{x}$ with a small angle $\theta$ with respect to 
$\widehat{\eta}$, stays close to it, and reaches $\widehat{y}$ with angle $\theta$ on the other 
side, crossing $\widehat{\eta}$ just once at some point in the lower-half plane (see left side of 
Figure \ref{pic_nonisobub2}). 
This arc can be chosen to sit inside a fundamental domain for $\rho(\gamma)$, so that it gives an 
embedded arc on $Gr(\sigma_\rho,\gamma)$ starting at $x^-$, 
reaching $\gamma^+$ at a point $z^+$ close to $x^+$ and ending at $y$. Changing the value of 
$\theta$ in some small interval $]-\varepsilon,\varepsilon[$ we obtain a family of embedded arcs 
$\alpha_\theta$ in $Gr(\sigma_\rho,\gamma)$ which are isotopic relative to the endpoints $x^-,y$ 
and are all injectively developed, except $\alpha_0=\eta$. 

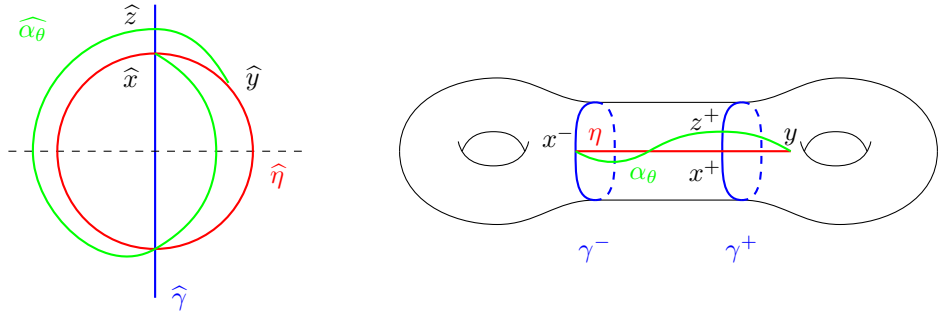
\begin{figure}[h]

\begin{center}
\begin{tikzpicture}[scale=0.65]

\draw[blue,thick] (5,3) to[out=-90,in=90] (5,-3);
\draw[dashed] (2,0) to[out=0,in=180] (8,0);
\draw[red,thick] (5,2) to[out=0,in=90] (7,0) to[out=-90,in=0] (5,-2) to[out=180,in=-90] 
(3,0) to[out=90,in=180] (5,2);
\draw[green,thick] (5,2) to[out=-30,in=90] (6.25,0) to[out=-90,in=30] (5,-2) to[out=210,in=-90] 
(2.5,0) to[out=90,in=180] (5,2.5) to[out=0,in=120] (6.5,1.4);
\node at (5.5,-3) {\color{blue} $\widehat{\gamma}$};
\node at (7.5,-0.5) {\color{red} $\widehat{\eta}$};
\node at (2.5,2.5) {\color{green} $\widehat{\alpha_\theta}$};
\node at (4.5,1.5) {$\widehat{x}$};
\node at (7,1.5) {$\widehat{y}$};
\node at (4.5,2.8) {$\widehat{z}$};

\draw[xshift=10cm] (0,0) to[out=90,in=180] (2,1.5) to[out=0,in=180] (4,1) to[out=0,in=180] (7,1) 
to[out=0,in=180] (9,1.5) to[out=0,in=90] (11,0) ;
\draw[xshift=10cm,yscale=-1] (0,0) (0,0) to[out=90,in=180] (2,1.5) to[out=0,in=180] (4,1) 
to[out=0,in=180] (7,1) to[out=0,in=180] (9,1.5) to[out=0,in=90] (11,0) ;
\draw[xshift=10cm,xscale=0.8]  (1.6,0) to[out=65,in=180] (2.4,0.4) to[out=0,in=115] (3.2,0);
\draw[xshift=10cm,xscale=0.8]  (1.5,0.2) to[out=-75,in=180] (2.4,-0.3) to[out=0,in=-105] (3.3,0.2);
\draw[xshift=17cm,xscale=0.8]  (1.6,0) to[out=65,in=180] (2.4,0.4) to[out=0,in=115] (3.2,0);
\draw[xshift=17cm,xscale=0.8]  (1.5,0.2) to[out=-75,in=180] (2.4,-0.3) to[out=0,in=-105] (3.3,0.2);

\draw[xshift=10cm,blue,thick] (4,1) to[out=180,in=90] (3.6,0) to[out=-90,in=180] (4,-1);
\draw[xshift=10cm,blue,thick,dashed] (4,1) to[out=0,in=90] (4.4,0) to[out=-90,in=0] (4,-1);
\draw[xshift=13cm,blue,thick] (4,1) to[out=180,in=90] (3.6,0) to[out=-90,in=180] (4,-1);
\draw[xshift=13cm,blue,thick,dashed] (4,1) to[out=0,in=90] (4.4,0) to[out=-90,in=0] (4,-1);

\node at (14,-2) {\color{blue} $\gamma^-$};
\node at (17,-2) {\color{blue} $\gamma^+$};


\draw[red,thick] (13.6,0) to[out=0,in=180] (18,0);
\node at (13.25,0.3) { $x^-$};
 \node at (16.25,-0.4) { $x^+$};
 \node at (16.25,0.65) { $z^+$};
\node at (18,0.3) { $y$};
\node at (14,0.3) {\color{red} $\eta$};
\draw[green,thick](13.6,0) to[out=-30,in=210] (15.1,0) to[out=30,in=180] (16.6,0.4) 
to[out=0,in=150] (18,0);
\node at (15,-0.5) {\color{green} $\alpha_\theta$};

\end{tikzpicture}

\end{center}
\caption{The bubbleable arc $\alpha_\theta$ in $\cp$ and 
$Gr(\sigma_\rho,\gamma)$}\label{pic_nonisobub2}
\end{figure}

Fix now some small $\theta $ and consider the BPS $\sigma_\pm = 
Bub(Gr(\sigma_\rho,\gamma),\alpha_{\pm \theta})$ obtained by bubbling along $\alpha_{\pm \theta}$. 
We now proceed to show that these two BPSs are not isomorphic: they can be distinguished by looking 
at the configuration of certain curves, which we now define.
The first curve we need is the analytic continuation of $\gamma^+$:  
we extend it inside the bubble by following its developed image. The result is a curve which  
reaches $x^-$, and we still denote it by $\gamma^+$. 
To define the other curve, let us recall from \cite[\S 3]{CDF} that a BPS with Fuchsian holonomy 
canonically decomposes into subsurfaces endowed with (possibly branched) complete hyperbolic 
metrics (see also Lemma \ref{geometryinqfholonomy} below for more details). Then the curve we need is the 
unique geodesic $\delta$ between $x^-$ and $y$ with respect to this metric, which develops 
isometrically onto 
the developed image of the original geodesic segment $\eta$ of $\sigma_\rho$. Notice that the whole 
construction can be made in such a way that this is indeed the shortest geodesic between its 
endpoints, just by taking the segment $\eta$ on $\sigma_\rho$ to be suitably shorter than the 
systole of $\sigma_\rho$.
Now we look at the tangent space at $x^-$. The tangent vector to $\gamma^+$ at 
$x^-$ sits on the right or on the left of the tangent vector to $\delta$ (with respect to the 
underlying orientation of $S$) depending on the fact that we look at $\sigma_+$ or at $\sigma_-$.
But any projective isomorphism between the two structures should be in particular orientation 
preserving at $x^-$.
\end{example}

\section{Geometric decomposition in \qf holonomy}\label{s_geometricdecomposition}

We now restrict our attention to structures whose holonomy preserves a decomposition of the model 
space $\cp$ into two disks separated by a Jordan curve. As observed in \cite{GO} and \cite{CDF}, 
the key feature of structures with such a representation is the presence of a canonical 
decomposition of the surface into subsurfaces which carry complete (possibly branched) hyperbolic 
structures with ideal boundary. The purpose of this section is to give a description of the 
components that can appear in such a decomposition, in the spirit of Goldman's work in \cite{GO}. 
Let us begin by recalling some definitions and known constructions.
\begin{definition}\label{def_qfgrp}
A \textbf{Fuchsian} (respectively \textbf{quasi-Fuchsian}) \textbf{group} is a subgroup of $\pslc$ 
 whose limit set in $\cp$ is $\rp$ (respectively a Jordan curve).
 \end{definition}
In particular a finitely generated quasi-Fuchsian group $\Gamma$  preserves a decomposition $\cp=\Omega^+_\Gamma \cup \Lambda_\Gamma \cup 
\Omega^-_\Gamma$ of the Riemann sphere into a pair of disks $\Omega^\pm_\Gamma$ and a Jordan curve $\Lambda_\Gamma$, i.e. the two components of the domain of discontinuity and the limit set of $\Gamma$. When $\Gamma$ is Fuchsian this is just the decomposition $\cp=\mathcal{H}^+ \cup \rp \cup \mathcal{H}^+$, where $\mathcal{H}^\pm$ denote the upper and lower-half plane in $\cmp$.
 \begin{definition}\label{def_qfrep}
 A faithful representation $\rho:\pi_1(S)\hookrightarrow \pslc$ is a \textbf{Fuchsian} (respectively \textbf{quasi-Fuchsian}) \textbf{representation} if its image is a Fuchsian (respectively quasi-Fuchsian) subgroup and there exists an orientation preserving 
$\rho$-equivariant diffeomorphism $f:\widetilde{S}\to \Omega^+_{\rho(\pi_1(S))}$. A structure $\sigma \in \mathcal{BP}(S)$ is said to be Fuchsian or quasi-Fuchsian when its holonomy is.
\end{definition}
By classical results (see for instance \cite[Theorem 4]{BE}) finitely generated quasi-Fuchsian surface groups are obtained as quasi-conformal deformations of Fuchsian ones; more precisely, given a quasi-Fuchsian representation $\rho$, there exist a Fuchsian representation $\rho_0$ and an orientation preserving quasi-conformal self-homeomorphism of $\cp$ that conjugates $\rho$ to $\rho_0$. This allows to extend the main results established in the following to quasi-Fuchsian representations, but we restrict to the Fuchsian case for the sake of simplicity.\par

Given a Fuchsian representation $\rho$, by definition we have an orientation preserving 
$\rho$-equivariant diffeomorphism $f:\widetilde{S}\to \mathcal{H}^+$. This descends to an 
orientation preserving diffeomorphism $F:S\to \mathcal{H}^+ / Im(\rho)$, giving us a (marked) 
unbranched complete hyperbolic structure on $S$ with holonomy $\rho$; we can use it as a base point 
in the moduli space $\mathcal{M}_\rho$, so we give it a special name (compare \cite[Definition 1.5]{GO} and \cite[Definition 2.3]{CDF}).
\begin{definition}
If  $\rho:\pi_1(S)\to \pslc$ is  a Fuchsian representation, then $\sigma_\rho=\mathcal{H}^+_\rho / Im(\rho)$ 
is called the \textbf{uniformizing structure} for $\rho$.
\end{definition}
More generally, if $dev:\widetilde{S}\to \cp$ is a developing map for a BPS on $S$  with \qf 
holonomy $\rho$, then the decomposition of the Riemann sphere induced by $\rho$ can be pulled back 
via the $dev$ to obtain a decomposition of $\widetilde{S}$. Since the developing map 
is ($\pi_1(S),\rho$)-equivariant, this decomposition is $\pi_1(S)$-invariant and thus descends to a 
decomposition of the surface into possibly disconnected subsurfaces $\sigma^+$ and 
$\sigma^-$ and a possibly disconnected curve $\sigma^\real$ defined as the subset of points 
developing to $\mathcal{H}^+$, $\mathcal{H}^-$ and $\rp$ respectively.
\begin{definition}\label{def_geomdecomp}
We will call $S=\sigma^+ \cup \sigma^\real \cup \sigma^-$ the \textbf{geometric decomposition} of 
$S$ with  respect to the BPS defined by the pair ($dev,\rho$); we will call $\sigma^\pm$ the 
positive/negative part of $S$ and $\sigma^\real$ the real curve of $S$. 
\end{definition}
We already observe at this point that, despite their apparent symmetry, the positive and negative 
part play a very different role in the geometry of $\sigma$, because of the special role played by 
$\mathcal{H}^+$ in the Definition \ref{def_qfrep} of Fuchsian 
representation. This phenomenon was already exploited by Goldman in the unbranched 
case (see \cite{GO}), and we will explore the branched case below. \par
Notice that a priori the decomposition of the surface depends not only on the representation, but  
also on the choice of a developing map. However this ambiguity can be fixed by choosing a 
representation $\rho$ in its conjugacy class, as explained in Remark \ref{uniquedev} above, since \qf 
representations are in particular non-elementary representations.
As a result, the decomposition of $S$ depends only on the structure $\sigma=\{(dev,\rho)\}$ and not 
on the choice of particular representatives. In particular many combinatorial properties of the 
geometric decomposition (such as the number and type of components, the adjacency pattern, 
the location of branch points,\dots) are well defined. The following was observed in \cite[\S 2]{GO} 
for the unbranched case and in \cite[\S 3]{CDF} for the branched case, and is the main feature of 
structures with \qf holonomy.
\begin{lemma}\label{geometryinqfholonomy}
If $S$ is endowed with a \qf BPS $\sigma$, then  $\sigma^\pm$ is a finite union of subsurfaces 
carrying complete hyperbolic metrics with cone points of angle $2\pi(k+1)$ corresponding to 
branch points of order $k$ of the BPS, and $\sigma^\real$ is a finite 1-dimensional CW-complex on 
$S$; moreover if branch points are not on the real curve, then $\sigma^\real$ is a finite union of 
simple closed curves with a $(\pslr,\rp)$-structure. 
\end{lemma}
Moreover this motivates the following terminology.
\begin{definition}
If $S$ is endowed with a \qf BPS $\sigma$, a connected component $C$ of $\sigma \setminus 
\sigma^\real$ will be called a geometric component of the decomposition; a connected component $C$ 
of $\sigma^\pm$ will be called a positive/negative component. A connected component of 
$\sigma^\real$ will be called a real component.
\end{definition}
Notice  that the components of the real curve can be canonically oriented by declaring that they 
have positive regions on the left and negative regions on the right. Some examples are in order.
\begin{example}
A hyperbolic structure on $S$ is an example of an unbranched projective structure  with Fuchsian 
holonomy. Any developing map is a diffeomorphism with the upper-half plane $\mathcal{H}^+$. The 
induced decomposition is $\sigma^+=S, \sigma^-=\varnothing=\sigma^\real$. Hence there is only one 
geometric component, which is the whole surface.
\end{example}
\begin{example}
If we graft a hyperbolic surface along a simple closed geodesic we obtain an example  of an 
unbranched projective structure with Fuchsian holonomy with surjective and non-injective 
developing map to $\cp$. There are a negative geometric annulus bounded by two essential simple 
closed real curves and two or one positive geometric components, depending on the fact that the 
geodesic we use is separating or not.
\end{example}
The main result in \cite{GO} claims that every unbranched structure with \qf holonomy arises via a 
multigrafting of the uniformizing hyperbolic structure; one of the key observations is the fact 
that geometric components of an unbranched structure can not be simply connected, i.e. they can 
not be disks. This completely fails for branched structures as the following easy example shows.
\begin{example}\label{ex_bubhyp}
If we bubble a hyperbolic surface along a simple arc we obtain an example of a branched projective 
structure with Fuchsian holonomy and with a negative geometric disk bounded by a contractible 
simple closed real curve and one positive geometric component containing the two branch points.
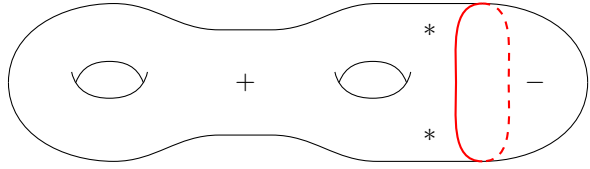
\begin{figure}[h]

\begin{center}
\begin{tikzpicture}[scale=0.7]

\draw (2,0) to[out=90,in=180] (4,1.5) to[out=0,in=180] (6,1) to[out=0,in=180] (7,1) 
to[out=0,in=180] (9,1.5) to[out=0,in=180] (11,1.5) to[out=0,in=90] (13,0) ;
\draw[yscale=-1] (2,0) to[out=90,in=180] (4,1.5) to[out=0,in=180] (6,1) 
to[out=0,in=180] (7,1) to[out=0,in=180] (9,1.5) to[out=0,in=180] (11,1.5) to[out=0,in=90] (13,0) ;
\draw[xshift=2cm,xscale=0.8]  (1.6,0) to[out=65,in=180] (2.4,0.4) to[out=0,in=115] (3.2,0);
\draw[xshift=2cm,xscale=0.8]  (1.5,0.2) to[out=-75,in=180] (2.4,-0.3) to[out=0,in=-105] (3.3,0.2);
\draw[xshift=7cm,xscale=0.8]  (1.6,0) to[out=65,in=180] (2.4,0.4) to[out=0,in=115] (3.2,0);
\draw[xshift=7cm,xscale=0.8]  (1.5,0.2) to[out=-75,in=180] (2.4,-0.3) to[out=0,in=-105] (3.3,0.2);

\draw[red,thick] (11,1.5) to[out=180,in=90] (10.5,0) to[out=-90,in=180] (11,-1.5);
\draw[red,thick,dashed] (11,1.5) to[out=0,in=90] (11.5,0) to[out=-90,in=0] (11,-1.5);

\node at (6.5,0) {$+$};
\node at (12,0) {$-$};
\node at (10,-1) {$*$};
\node at (10,+1) {$*$};

\end{tikzpicture}

\end{center}
\caption{Geometric decomposition of a bubbling on $\sigma_\rho$.}\label{pic_geomdecompbub}
\end{figure}
\end{example}
Our purpose here is to show that also for branched structures it is possible to obtain a control of 
the behaviour of the negative components. The location of branch points with respect to the 
geometric decomposition is of course something we 
want to care about in the following, therefore we introduce the following definitions.
\begin{definition}\label{def_geomrealbranch}
 Let $\sigma$ be a \qf BPS. A branch point of $\sigma$ is said to be geometric (respectively 
real) if it belongs to $\sigma^\pm$ (respectively to $\sigma^\real$). The structure is said to be 
\textbf{geometrically branched} if all its branch points are geometric.

We will denote by $\mathcal{M}_{k,\rho}^\real$ the subspace of non-geometrically branched structures of $\mathcal{M}_{k,\rho}$, i.e. those with at least one real branch point.
\end{definition}
Notice that up to a very small movement of branch points, we can always assume that the branch 
points do not belong to the real curve $\sigma^\real$; more precisely, $\mathcal{M}_{k,\rho}^\real$ 
has real codimension $1$ inside the $k$-dimensional complex manifold $\mathcal{M}_{k,\rho}$. From 
now on we focus on geometrically branched structures; for these ones some index formulae are 
available, which link the geometry and the topology of the components of the 
geometric decomposition. We recall here the needed terminology (see \cite[\S 3-4]{CDF} for more 
details.)
\begin{definition}
 Let $\sigma$ be a geometrically branched BPS and $l$ be a real component on it. Let $p \in 
\rp$ be a fix point of $\rho(l)$ and $\widetilde{l}$ is any lift of $l$. 
The \textbf{index} of the induced real projective structure on $l$ is the integer 
$I(l)=\#\left( \{ dev^{-1}_{|_{\widetilde{l}}}(p) \} / <l> \right)$.
\end{definition}
The index of a real component can be thought as a degree of the restriction of the developing 
map to it, as a map with values in the limit set of $\rho$, and it can a priori assume any value. 
However if $\rho(l)$ is trivial then the index must be strictly positive: this follows by the 
classification of $\rp$-structures on $S^1$ given in \cite[Proposition 3.2]{CDF}, which we recall 
for future reference.
\begin{lemma}\label{trivialholopositiveindex}
 Two unbranched $\rp$-structures on an oriented circle with non-elliptic holonomy are isomorphic if 
and only if they have the same index and holonomy. The only case which does not occur is 
the case of index $0$ and trivial holonomy.
\end{lemma}
\begin{definition}
For a \qf representation $\rho$ let $E_\rho$ be the induced flat $\rp$-bundle on $S$. For any 
subsurface $i:C\hookrightarrow S$ we denote by $\rho_C$ the restriction of $\rho$ to $i_*\pi_1(C)$. 
For any component $l\subset \partial C$ we define a section $s_\rho:l \to 
{E_\rho}{|_{l}}$ by choosing the flat section passing through a fixed point of 
$\rho(l)$. Then the \textbf{Euler class} $eu$ of $\rho_C$ is defined to be the Euler class of 
the bundle $E_{\rho_C}={E_\rho}{|_{C}}$ with respect to this choice of boundary sections.
\end{definition}
Finally we say that a subsurface $C\subset S$ is incompressible if the inclusion is injective on 
fundamental groups or, equivalently, if all the boundary curves are essential (i.e. not 
nullhomotopic) in $S$. The following index formulae hold.
\begin{theorem}(\cite[Theorem 4.1-5]{CDF})\label{indexformula}
 Let $\sigma \in \mathcal{M}_\rho$ be geometrically branched. Let $C\subset \sigma^\pm$ be a 
geometric component containing $k_C$ branch points (counted with multiplicity) and with $\partial C 
=\{l_1,\dots,l_n\}\subset \sigma^\real$. Then
$$\pm eu(\rho_C)=\chi(C)+k_C-\sum_{i=1}^n I(l_i)$$
Moreover if $C$ is incompressible (e.g. $C=S$) then $eu(\rho_C)=\chi(C)$.
\end{theorem}
Under the same hypothesis of Theorem \ref{indexformula} the following can be deduced
\begin{corollary}
 If $k^\pm$ denotes the number of positive/negative branch points of $\sigma$, then 
$2\chi(\sigma^-)=k^+-k^-$.
\end{corollary}
In particular in \qf holonomy there is always an even number of branch points, so that 
$\mathcal{M}_{2k+1,\rho}$ are all empty. 
\begin{example}\label{ex_bubhaveposindex}
 If $\sigma$ is unbranched, every geometric component carries an unbranched complete hyperbolic 
metric by Lemma \ref{geometryinqfholonomy}; as a consequence all real curves have index $0$, and in 
particular they are essential by Lemma \ref{trivialholopositiveindex}. On the other hand in the branched 
case real curves can have positive index and be non-essential; as an example consider the simple 
bubbling of Example \ref{ex_bubhyp}, where there is exactly one contractible real component with index $1$.
\end{example}

\subsection{Locating branch points}
We have observed in Example \ref{ex_bubhaveposindex} that genuinely branched structures can have real 
curves of positive index. Roughly speaking, if this occurs then branch points must live in the 
geometric components adjacent to the real curves of positive index. This section aims at making 
this statement more precise. 
We begin by noticing that even if the structure has branch points, nevertheless unbranched 
components are quite well behaved. 
\begin{lemma}\label{unbrcpt}
Let $\sigma \in \mathcal{M}_{k,\rho}$ be geometrically branched.
If $C\subset \sigma^\pm$ is an unbranched component then either it is a disk or it is 
incompressible. 
Moreover if it is negative and incompressible, then it is an incompressible annulus.
\end{lemma}
\begin{proof} We already know that unbranched disks can occur. If $C$ is not a disk and is not branched, 
then it 
carries a complete hyperbolic structure such that the index of each boundary component is zero; by Lemma
\ref{trivialholopositiveindex} we know that it can not have trivial holonomy. But \qf 
representations are in particular injective, hence this implies that each boundary component must 
be essential in the surface $S$, hence $C$ is incompressible. So we can apply the index formula and 
obtain $\chi(C)=eu(\rho_C)=\pm \left( \chi(C)+k_C-\sum_{l \subset \partial C} I(l) 
\right)=\pm \chi(C)$, where the sign depends on the sign of $C$. In the case $C$ is negative this 
implies $\chi(C)=0$.
\end{proof}
This is the first manifestation of the asymmetry between positive and negative 
regions hinted at before, and which is a consequence of the special role played by $\mathcal H^+$ 
in the Definition \ref{def_qfrep} of \qf representation. Notice that we get a useless identity in 
the case of a positive component. The following easy observation provides a first step to locate 
branch points with respect to $\sigma^\real$.
\begin{lemma}\label{posindex}
Let $\sigma \in \mathcal{M}_{k,\rho}$ be geometrically branched.
Let $l$ be a real component and $C,C'$ be the components of $\sigma^\pm$ which are adjacent 
along $l$. If $I(l)\geq 1$ then
\begin{enumerate}
\item at most one of $C,C'$ is a disk;
\item any non disk component is branched;
\item at least one of $C,C'$ is branched.
\end{enumerate}
\end{lemma}
\begin{proof}
To prove (1) observe that $C,C'$ can not both be disks, otherwise we would get an embedded sphere 
in $S$, which is always assumed to have genus $g\geq 2$. To get (2) observe that an unbranched 
component with ideal boundary with positive index is  necessarily a disk isometric to 
$\hyp^2$: indeed such a component carries a complete hyperbolic structure, hence is a quotient of 
$\hyp^2$ by some group $\Gamma$ and  as soon as $\Gamma\neq id$ we see that the index of the real 
boundaries is $0$. Finally (3) follows from (1) and (2).
\end{proof}
\begin{corollary}\label{adjdisk}
Let $\sigma \in \mathcal{M}_{k,\rho}$ be geometrically branched.
The (unique) component adjacent to a disk of the geometric decomposition of $\sigma$ is branched.
\end{corollary}
\begin{proof}
The boundary of a disk has always strictly positive index. Since $S$ is not a sphere, the adjacent 
component can not be a disk, therefore it is branched.
\end{proof}
Notice however that it may happen that none of the two components adjacent to a real component with 
positive index is a disk.
\begin{example}\label{ex1+}
Let $\gamma$ be a simple closed geodesic on the uniformizing structure $\sigma_\rho$, and let 
$\sigma'=Gr(\sigma_\rho,\gamma)$. Then pick a bubbleable arc $\beta \subset \sigma'$ which 
intersects exactly once the real curve of $\sigma'$ and let $\sigma''=Bub(\sigma',\beta)$. Then 
$\sigma''$ has an essential real component of index $1$ such that both adjacent components are non 
disks (both have non-positive Euler characteristic) and are branched. See Figure 
\ref{pic_exlocatebr}, left side.
\end{example}
On the other hand there are structures with negative components with essential boundary with index 
$0$ which are nevertheless branched.
\begin{example}\label{ex2-}
Let $\gamma$ be a non-separating simple closed geodesic on the uniformizing structure 
$\sigma_\rho$, 
and let $\sigma'=Gr(\sigma_\rho,\gamma)$. Then let $\beta \subset \sigma'$ be a bubbleable arc with 
endpoints inside the negative annulus but which is not itself contained inside the 
negative annulus and let $\sigma''=Bub(\sigma',\beta)$. Then $\sigma''$ has one real component of 
index $0$, a positive unbranched incompressible component and a negative incompressible component 
of Euler characteristic -1 containing both branch points. See Figure \ref{pic_exlocatebr}, right 
side.
\end{example}
\begin{figure}[h]

\begin{center}
\begin{tikzpicture}[scale=0.6,xscale=-1]
\draw (0,0) to[out=90,in=180] (2,1.5) to[out=0,in=180] (3,1.5) to[out=0,in=180] (5,1) 
to[out=0,in=180] (7,1.5) to[out=0,in=90] (9,0) ;
\draw[yscale=-1] (0,0) to[out=90,in=180] (2,1.5) to[out=0,in=180] (3,1.5) to[out=0,in=180] (5,1) 
to[out=0,in=180] (7,1.5) to[out=0,in=90] (9,0) ;
\draw[xscale=0.8]  (1.6,0) to[out=65,in=180] (2.4,0.4) to[out=0,in=180] (3.4,0.4) 
to[out=0,in=115] (4.2,0);
\draw[xscale=0.8]  (1.5,0.2) to[out=-75,in=180] (2.4,-0.3) to[out=0,in=180] (3.4,-0.3) 
to[out=0,in=-105] (4.3,0.2);
\draw[xshift=5cm,xscale=0.8]  (1.6,0) to[out=65,in=180] (2.4,0.4) to[out=0,in=115] (3.2,0);
\draw[xshift=5cm,xscale=0.8]  (1.5,0.2) to[out=-75,in=180] (2.4,-0.3) to[out=0,in=-105] (3.3,0.2);

\draw[red,thick] (1.8,-0.3) to[out=180,in=90] (1.6,-0.9) to[out=-90,in=180] (1.8,-1.5);
\draw[red,dashed] (1.8,-0.3) to[out=0,in=90] (2,-0.9) to[out=-90,in=0] (1.8,-1.5);
\draw[red,thick,xshift=1cm] (1.8,-0.3) to[out=180,in=90] (1.6,-0.9) to[out=-90,in=180] (1.8,-1.5);
\draw[red,dashed,xshift=1cm] (1.8,-0.3) to[out=0,in=90] (2,-0.9) to[out=-90,in=0] (1.8,-1.5);

\draw[blue,thick] (1.8,-0.9) to[out=180,in=-90] (0.75,0) to[out=90,in=180] (2.3,0.9) 
to[out=0,in=90] (4,0) to[out=-90,in=0] (2.4,-0.9);

\draw[xshift=10cm] (0,0) to[out=90,in=180] (2,1.5) to[out=0,in=180] (4,1) to[out=0,in=180] (7,1) 
to[out=0,in=180] (9,1.5) to[out=0,in=90] (11,0) ;
\draw[xshift=10cm,yscale=-1] (0,0) (0,0) to[out=90,in=180] (2,1.5) to[out=0,in=180] (4,1) 
to[out=0,in=180] (7,1) to[out=0,in=180] (9,1.5) to[out=0,in=90] (11,0) ;
\draw[xshift=10cm,xscale=0.8]  (1.6,0) to[out=65,in=180] (2.4,0.4) to[out=0,in=115] (3.2,0);
\draw[xshift=10cm,xscale=0.8]  (1.5,0.2) to[out=-75,in=180] (2.4,-0.3) to[out=0,in=-105] (3.3,0.2);
\draw[xshift=17cm,xscale=0.8]  (1.6,0) to[out=65,in=180] (2.4,0.4) to[out=0,in=115] (3.2,0);
\draw[xshift=17cm,xscale=0.8]  (1.5,0.2) to[out=-75,in=180] (2.4,-0.3) to[out=0,in=-105] (3.3,0.2);

\draw[xshift=10cm,red,thick] (4,1) to[out=180,in=90] (3.6,0) to[out=-90,in=180] (4,-1);
\draw[xshift=10cm,red,thick,dashed] (4,1) to[out=0,in=90] (4.4,0) to[out=-90,in=0] (4,-1);
\draw[xshift=13cm,red,thick] (4,1) to[out=180,in=90] (3.6,0) to[out=-90,in=180] (4,-1);
\draw[xshift=13cm,red,thick,dashed] (4,1) to[out=0,in=90] (4.4,0) to[out=-90,in=0] (4,-1);

\draw[blue,thick] (15.5,0) to[out=0,in=180] (18,0);

\end{tikzpicture}

\end{center}
  \caption{Figures for Examples \ref{ex1+} and \ref{ex2-}.}\label{pic_exlocatebr}
  \end{figure}

From the above results we obtain in particular a bound on the number of branch points contained 
inside a disk of the geometric decomposition of a \qf structure $\sigma \in \mathcal{M}_{k,\rho}$.
\begin{proposition}\label{brdisk}
Let $\sigma \in \mathcal{M}_{k,\rho}$ be geometrically branched and $k\geq 2$.
If a geometric component $D\subset \sigma^\pm$ is a disk of branching order $k_D$, then $k_D \leq 
k-2$.
\end{proposition}
\begin{proof}
By Corollary \ref{adjdisk} we already know that a disk can not contain all the branching. So we assume by 
contradiction that it has branching order $k_D=k-1$. Since $D$ is contractible, by Theorem
\ref{indexformula} its boundary  is a real component $l$ of 
index $I(l)=k_D+1=k$. Let $C$ be the component adjacent to $D$; then we know it is branched by Corollary
\ref{adjdisk}, so $k_C\geq 1$. Indeed $k_D=k-1 $ implies that $ k_C=1$. The boundary of $C$ a 
priori could contain also $m$ more non-essential boundary components and $n$ essential ones. Notice 
that all components of $\sigma^\pm$ different from $C,D$ are unbranched, simply because $C\cup D$ 
contains all the branching. \par
Therefore if $l'\neq l$ is a non-essential component of $\partial C$, then $\rho(l')=id$ hence 
$I(l')\geq1$ by Lemma \ref{trivialholopositiveindex}, and then by Lemma \ref{posindex} the geometric component 
after it must be an unbranched disk $D'$ and $l'$ must have index $1$. Let $l,l'_1,\dots, l'_m$ be 
the non-essential components of $\partial C$, and $D,D'_1,\dots,D'_m$ the corresponding disks; then 
$I(l)=k$ but $I(l'_i)=1,k_{D'_i}=0$ for 
$i=1,\dots,m$. 
On the other hand, if $l''$ is an essential boundary component, then the geometric 
component after it is a non-simply connected complete hyperbolic surface, hence $I(l'')=0$.\par
Now observe that the subsurface $E=C\cup D \cup D'_1 \cup \dots \cup D'_m$ has essential boundary 
by construction, hence it is incompressible. Therefore the index formula (Theorem \ref{indexformula}) gives 
us that
$$eu(\rho_E)=\chi(E) = \chi (C) + \chi(D)+\sum_{i=1}^m \chi(D'_i)= \chi (C) 
+ 1+m $$
On the other hand, we obtain, again by Theorem
\ref{indexformula} and the fact that disks have trivial Euler class, that
$$eu(\rho_E)= eu(\rho_C) + 
eu(\rho_D)+\sum_{i=1}^m eu(\rho_{D'_i})= 
eu(\rho_C) = $$
$$= \pm \left(\chi(C)+k_C-I(l)-\sum_{i=1}^m I(l'_i)-\sum_{j=1}^m I(l''_j)\right)= \pm 
(\chi(C) +1 -k-m)  $$
where the sign depends on the sign of $C$ (hence of that of $D$). We are now going to compare the 
two expressions for the Euler class of $E$. If $C\subset \sigma^+$ then we 
get $2m+k=0$ which is absurd since $m\geq 0, k\geq 2$. If $C\subset \sigma^-$ then we get 
$2\chi(C)=k-2\geq 0$. But $C$ can not be a disk, hence $\chi(C)=0$, i.e. $C$ is an annulus. Its 
boundary consists of $l$ and another curve $l'$ homotopic to it; so $l'$ is non-essential too, 
hence of positive index. The component adjacent to $l'$ can not be a disk, otherwise $S$ would have 
genus $g=0$, hence it must be branched; but by construction all branch points live in $C\cup D$, so 
we have a contradiction.
\end{proof}

Notice that so far $C$ could be either positive or negative. Indeed, by performing suitable 
bubbling, we can find structures with either positive or negative disks, either branched or not. 
We recall the following useful lemma, which was proved in \cite[Lemma 10.3]{CDF} for the positive 
part; here we just show that the same proof provides an interesting equality for the negative part 
too.
\begin{lemma}\label{alltogether}
Let $\sigma \in \mathcal{M}_{k,\rho}$ be geometrically branched.
If all branch points live in $\sigma^+$ and $C\subset \sigma^+$ is a branched component with $n$ 
adjacent disks, then $k_C=2n$. If all branch points live in $\sigma^-$ and $C\subset \sigma^-$ is a 
branched component then $k_C=-2\chi(C)$.
\end{lemma}
\begin{proof}
Suppose all branch points live in the positive part or in the negative part, and let $C$ be a 
branched component. The hypothesis implies that all components adjacent to $C$ are unbranched, 
therefore by Lemma \ref{unbrcpt} we have the following dichotomy for a real curve in the boundary of $C$: 
either it has index $0$ and is essential, or it has index $1$ and  bounds a disk. Let 
$l_1,\dots,l_n$ be the non-essential boundary components of $C$ and let $D_1,\dots,D_n$ be the 
adjacent disks. The subsurface $E=C \cup D_1 \cup \dots \cup D_n$ is clearly 
incompressible. By Theorem \ref{indexformula} and the fact that disks have trivial Euler class we obtain
$$\chi(C)+n=\chi(E)=eu(\rho_E)=eu(\rho_C)= eu(\rho_C)+\sum_{i=1}^n eu(\rho_{D_i})=$$
$$\pm \left( \chi(C)+k_C -\sum_{i=1}^n I(l_i)  \right) = \pm \left( \chi(C)+k_C -n \right)  $$
from which the statement follows.
\end{proof}

\subsection{Classification of components for BPSs with $k=2$ branch points}
When we have only two branch points, we can obtain a strong control on the behaviour of real curves 
of positive index. This can be used to obtain a classification of the components that can appear in 
the geometric decomposition of a structure. As before we assume branch points are not on the real 
curve, so that the index formulae Theorem \ref{indexformula} can be used.\par
In Lemma \ref{unbrcpt} we observed that in general an unbranched negative component which is not a 
disk is automatically an incompressible annulus. For structures with two branch points we can 
obtain a precise statement also about branched negative incompressible components.
\begin{lemma}\label{negbranch}
Let $\sigma \in \mathcal{M}_{2,\rho}$ be geometrically branched.
 Let $C\subset \sigma^-$ be a branched negative incompressible component containing $k_C$ branch 
points. Then 
 \begin{enumerate}
  \item either $k_C=1$, $C$ is an annulus with $\partial C = l \cup l'$ such that $I(l)=0,I(l)=1$
  \item or $k_C=2$, $C$ is a pant or a once-holed torus and $\forall \ l \subset 
\partial C$ we have $I(l)=0$
 \end{enumerate}
\end{lemma}
\begin{proof} 
Since $C$ is incompressible we can applying the index formula and we get
$$-\chi(C)=-eu(\rho_C)=\chi(C)+k_C-\sum_{l\subset \partial C} I(l) 
\Rightarrow 2\chi(C)+k_C=\sum I(\gamma)\geq 0$$
and here we look for integer solutions with the constraints that $\chi(C)\leq 0$ (being 
incompressible, $C$ is not a disk) and $k_C\leq 2$. 
We see that the only possibilities are the following
\begin{enumerate}
 \item  $k_C=0, \chi(C)=0$, so that $C$ is an unbranched annulus (which we discard, since $C$ is 
assumed to be branched)
 \item  $k_C=1, \chi(C)=0$, so that $C$ is an annulus; we get $\sum I(l)=1$, which means that one 
boundary component has index $0$ and the other has index $1$
\item $k_C=2,\chi(C)=0$, so that $C$ is again an annulus and $\sum I(l)=2$; in particular there is 
a boundary with positive index and the adjacent component should be branched, but $C$ already 
contains all the branching (so we do not have this possibility)
\item $k_C=2,\chi(C)=-1$, and we have $\sum I(l)=0$, which implies that all boundaries have zero 
index.
\end{enumerate}
\end{proof}
To do a similar study for positive branched components we need some preliminary results. A 
straightforward consequence of Proposition \ref{brdisk} is that disks are always unbranched when we have only 
two branch points; in particular a real component bounding a geometric disk has index $1$. We want 
to prove an analogous statement for essential real components.
\begin{lemma}\label{badcpt}
Let $\sigma \in \mathcal{M}_{2,\rho}$ be geometrically branched.
If a component $C\subset \sigma^\pm$ 
is not a disk and contains a single simple branch point, 
then the inclusion $i:C\hookrightarrow S$ can not be nullhomotopic (i.e. $i_*(\pi_1(C))\subset 
\pi_1(S)$ can not be the trivial subgroup).
\end{lemma}
\begin{proof}
By contradiction assume $i_*(\pi_1(C))\subset \pi_1(S)$ is trivial. In particular $C$ must have 
genus $0$ and its boundary must consist of $m\geq 2$ (it is not a disk)  non-essential boundary 
components $l_1,\dots,l_m$ with index $I(l_i)\geq 1$. Since $i_*(\pi_1(C))$ is trivial in 
$\pi_1(S)$, the flat bundle associated to $\rho$ is trivial on $C$, hence the Euler class vanishes. 
Applying the index formula we obtain
$$0=\pm eu(\rho_C)=\chi(C)+k_C-\sum_{i=1}^m I(l_i)\leq 2-m+k_C-m \leq 
k_C-2$$
which contrasts with the fact that $k_C=1$.
\end{proof}

\begin{proposition}\label{indexbound}
Let $\sigma \in \mathcal{M}_{2,\rho}$ be geometrically branched.
If $l \subset \sigma^\real$ is any real component, then $I(l)\leq 1$.
\end{proposition}
\begin{proof}
Suppose by contradiction we have a real curve $l_0 \subset S_\real$ of index $I(l_0)\geq 2$.
We distinguish two cases.\par
In the case $l_0$ is homotopically trivial, it bounds exactly one subsurface $D$ homeomorphic to a 
disk one one side and another subsurface $S'$ which is not a disk on the other side. This 
subsurface $D$ can either be a geometric disk, or it can consist of more than just one single 
geometric component. In the first case it is unbranched by Proposition \ref{brdisk} hence $l_0$ should have 
index 1; in 
the second case the geometric component $C$ of $D$ which has $l_0$ in its boundary is a non disk 
component, hence it must be branched; since $S'$ must be branched as well by Lemma \ref{posindex}, $C$ 
contains exactly one branch point, but then Lemma \ref{badcpt} applies and we get a contradiction with 
the 
fact that $C$ is contained in a disk (i.e. with the fact that its inclusion is homotopically 
trivial).\par
For the second case, suppose  $l$ is essential. Let us call $C^\pm$ the adjacent geometric 
components. Then $C^\pm$ are branched by Lemma \ref{posindex}; more precisely $k_{C^\pm}=1$, they are not 
disks since $l_0$ is essential and all other components are unbranched, since $C^+ \cup C^-$ 
contains all the branching. The two components $C^\pm$ may have $m\geq 0$ more boundaries in 
common, let us call them $l_1,\dots,l_m$. Moreover each of them can have more boundary components, 
either essential or not. Let us focus on $C^+$; its boundary consists of $l_0, l_1,\dots,l_m$ and 
possibly of some other non-essential components $l'_1,\dots,l'_n$ and some  essential ones 
$l''_1,\dots,l''_p$, for some $n,p\geq 0$. Once again, the non-essential components 
$l'_1,\dots,l'_n$ must bound unbranched disks $D'_1,\dots,D'_n$ (hence they have index $1$), and 
the 
essential components $l''_1,\dots,l''_p$ must bound unbranched components which are not disks 
(hence 
they have index $0$ and are essential).\\
We consider the subsurface $E=C^+\cup D'_1 \dots D'_n$ and we see that it 
is incompressible: $l''_1,\dots,l''_p$ are essential by definition, $l_1,\dots,l_m$ are non-separating curves in $S$ ($C^+$ and $C^-$ are adjacent along $l_0$ in any case), hence they are 
essential as well, as soon as $m\geq 1$. The only case we need to check is when $m=0$, but we are 
currently discussing the case in which $l_0$ is essential.\par
Then we apply the index formula and get
$$eu(\rho_E)=\chi(E) = \chi (C^+) +\sum_{i=1}^n \chi(D'_i)= \chi (C) + n 
$$
On the other hand, as in the previous proofs, we obtain
$$eu(\rho_E)= eu(\rho_C) = 
\chi(C)+k_C-I(l_0)-\sum_{i=1}^m I(l_i)-\sum_{j=1}^n I(l'_j)-\sum_{h=1}^p 
I(l''_h)$$
$$  = \chi(C)+1 -I(l_0) -\sum_{i=1}^m I(l_i) - n $$
By comparing the two expressions we obtain that
$$2n + I(l_0) +\sum_{i=1}^m I(l_i) =1$$
Now we  have that the left hand side is a sum of non-negative integers and that $I(l_0)\geq 2$ 
by hypothesis, therefore in any case we reach an absurd.
\end{proof}

Now we can prove the following result about positive branched components, which is analogous to Lemma
\ref{negbranch} for the branched negative incompressible ones; a description of branched negative 
compressible components will follow from Theorem \ref{k=2} below.
\begin{lemma}\label{posbranch}
Let $\sigma \in \mathcal{M}_{2,\rho}$ be geometrically branched. Let $C\subset \sigma^+$ be a 
branched positive component. Then 
\begin{enumerate}
 \item if $C$ is incompressible then $k_C=1$ and there is a unique boundary curve of index $1$,
loxodromic holonomy and the component beyond it is branched;
 \item if $C$ is compressible then $k_C=2$ and there is a unique boundary curve of index $1$, 
trivial holonomy and the component beyond it is an unbranched disk.
 \end{enumerate}
\end{lemma}
\begin{proof}
If $C$ is incompressible then we apply the index formula and get $k_C = \sum_{l\subset \partial C} 
I(l)$. Moreover every boundary component is essential, and by Proposition \ref{indexbound} its index is at 
most $1$. Therefore we have exactly $k_C$ components of index $1$ (and possibly some components 
of index $0$). Being essential, they do not bound disks, hence the adjacent components are 
branched. In particular if $k_C=2$ then there are two boundaries with index $1$ and thus some 
branched component is adjacent to $C$; but $C$ already contains all the branching, hence $k_C=1$ 
and there is a unique real component of index $1$. Since it is essential and we are in \qf 
holonomy, the holonomy around the curve will be loxodromic. Of course the component beyond it 
is branched by Lemma \ref{posindex}.\\
If $C$ is compressible, then let us say there are $m\geq 1$ non-essential boundaries $l_1,\dots ,
l_m$ (which have index $1$ by Proposition \ref{indexbound}, since non-essential curves have strictly positive 
index by Lemma \ref{trivialholopositiveindex}) and $n\geq 0$ essential boundaries $l'_1,\dots ,l'_n$, 
$n_0$ of which have index $1$ (and 
the others have index $0$ by Proposition \ref{indexbound}).  Then we can cap $C$ with these adjacent negative 
disks and apply the index formula to the resulting incompressible subsurface $E$
$$\chi(C)+m=\chi(E)=eu( \rho_E)=eu( \rho_C)=\chi(C)+k_C-m-n_0$$
$$2m+n_0=k_C$$
Since $m\geq 1$ but $k_C\leq 2$, this implies that indeed $k_C=2, m=1,n_0=0$.
\end{proof}

The above study was focused on a single branched component, but now we go global with the help of Lemma
\ref{alltogether}.
\begin{theorem}\label{k=2}
 Let $S$ be a closed, connected and oriented surface of genus $g\geq 2$, $\rho:\pi_1(S)\to \pslc$ 
be a \qf representation and $\sigma\in \mathcal{M}_{2,\rho}$ be geometrically branched. Let $k^\pm$ 
denote the number of branch points in $\sigma^\pm$
\begin{enumerate}
\item If $k^+=2$ then both branch points live in the same positive component; more precisely 
there exists a unique negative unbranched disk and the branch points live in the positive 
component which is adjacent to it. 
\item If $k^-=2$ then both branch points live in the same negative component; more precisely 
there exists a negative component of Euler characteristic $-1$ containing both branch points. 
Moreover it has at most one non-essential boundary component (with trivial holonomy and index $1$), 
while all essential boundaries have loxodromic holonomy and index $0$.
\item If $k^+=k^-=1$ then the two branched components are adjacent along an essential real 
component with index $1$ and loxodromic holonomy; the negative branched component is an 
incompressible annulus.
\end{enumerate}
Moreover in each case all the other positive components are unbranched incompressible and all the 
other negative components are unbranched incompressible annuli and all the other real curves have 
index $0$.
\end{theorem}
\begin{proof}
We consider the three cases.
\begin{enumerate}
 \item We have $2\chi(\sigma^-)=k^+-k^-=2$, so $\chi(\sigma^-)=1$, thus there must be a negative 
disk $D$. 
Let $C$ be the positive component adjacent to $D$. By Lemma \ref{alltogether} $C$ contains $2$ (i.e. all) 
branch points and indeed there are no other negative disks.
\item We have $2\chi(\sigma^-)=k^+-k^-=-2$, so $\chi(\sigma^-)=-1$.  By Lemma \ref{alltogether} negative 
components 
are either unbranched incompressible annuli or components with Euler characteristic $-1$ and $2$ 
branch points; hence there is exactly one of the latter kind. If it is incompressible, then it has 
the required boundary behaviour by Lemma \ref{negbranch}. If it is a pair of pants and it has one non-essential boundary component, then the adjacent component is a disk (because it is unbranched), 
hence the index is $1$. If it had two non-essential boundaries, then also the third boundary would 
be non-essential, but then all the components adjacent to the three boundaries must be disk and $S$ 
would be a sphere, so this case is absurd.
\item Let $C$ be the positive branched component. Since it has only one branch point, by Lemma 
\ref{posbranch} it is incompressible and has a unique boundary component of index $1$ and 
hyperbolic holonomy. The negative component adjacent along it can not be a disk, hence it is 
branched, with one branch point. By Lemma \ref{negbranch} it is an incompressible annulus and the other 
boundary component has index $0$. Moreover notice that the only negative disks could appear at the 
boundary of $C$, but this is forbidden since it is incompressible.
\end{enumerate}
The rest of the statement follows  from the initial discussion: the non-branched components 
can not be disks, hence they are incompressible and with zero index boundary by Lemma \ref{unbrcpt}. The 
negative ones are annuli again by Lemma \ref{unbrcpt}. As a consequence all real curves have index 0, 
except the non-essential ones in the case $k^\pm=2$ and the curve separating the branch points in 
the case $k^+=1$, which have index $1$.
\end{proof}

This gives a description of negative branched components also in the compressible case, which was 
still missing so far.
\begin{rmk}\label{k=1tobbl}
A direct consequence of this classification is that in the case $k^+=1=k^-$ we can always satisfy 
the hypothesis of \cite[Theorem 7.1]{CDF}, hence we can move branch points without crossing the 
real 
curves to obtain a structure which is a bubbling over some unbranched structure. This is a key fact 
in the proof of the main theorem below (see Theorem \ref{mainbubblingthm}).
\end{rmk}
We conclude with the following minor but curious application of Theorem \ref{k=2}.
\begin{corollary}\label{k+s=odd}
Let $\sigma \in \mathcal{M}_{2,\rho}$ be geometrically branched. Then the number of branch 
points contained in $\sigma^+$ and the total number of real components always sum to an odd number.
\end{corollary}
\begin{proof}
If $k^+=0,2$ then there is a negative component of Euler characteristic $\pm 1$. In both cases it 
has an odd 
number of boundary components. All the other negative components are incompressible annuli. The 
total number of real components is therefore odd.
if $k^+=1$ then the positive branched component is incompressible and there is exactly one index 
$1$ real boundary, beyond which the negative branched component sits. And it is an annulus. All 
other negative components are annuli too, hence we have an even number of real components. 
\end{proof}

\section{BM-configurations}\label{s_BMconfig}
As observed above in Remark \ref{k=1tobbl}, when a structure with two branch points and \qf holonomy has a 
positive branch point and a negative one, then it can be slightly deformed inside 
$\mathcal{M}_{2,\rho}$ without changing the induced geometric decomposition so that a bubble 
appears. However it is not clear a priori whether this bubble is preserved when we keep deforming 
the structure to reach other regions of $\mathcal{M}_{2,\rho}$. In this section we study what 
happens when we try to move branch points along an embedded twin pair based at one of the vertices 
of a bubble (recall from Remark \ref{movingiscoord} than moving branch points provides local coordinates 
on $\mathcal{M}_{2,\rho}$). We find it useful to introduce the following notation: if 
$\mathcal{X}\subset \mathcal{M}_{2,\rho}$ then we denote by $\mathcal{BX}$ the subspace of 
$\mathcal{X}$ made of BPSs which are obtained via a bubbling over some unbranched structure from 
$\mathcal{M}_{0,\rho}$. At first we just pick a non-elementary representation $\rho$; we will 
specify when we will need to restrict to the \qf case.

\subsection{Standard BM-configurations}
We begin naively with the easy situation in which points can be moved without affecting the bubble.
\begin{definition}
 Let $\sigma \in \mathcal{BM}_{(1,1),\rho}$. A \textbf{BM-configuration} (Bubbling-Movement 
configuration) on $\sigma$ is the datum of a bubble $B$ together with an embedded twin pair $\mu$ 
based at a vertex $p$ of $B$. We denote the configuration by $(B,\mu,p)$.
\end{definition}
We introduce now the nicest type of BM-configuration, which will allow us to perform local 
deformations of the structure preserving the bubble.
\begin{definition}
 A BM-configuration $(B,\mu,p)$ on $\sigma \in \mathcal{BM}_{(1,1),\rho}$ is said to be a 
\textbf{standard BM-configuration} if either all the arcs are disjoint and disjointly developed 
outside the obvious intersections (i.e. $\partial B \cap \mu = \{p\}$ and $dev(\partial B)\cap 
dev(\mu)=\{dev(p)\}$) or the embedded twin pair is entirely contained 
in the boundary of the bubble (i.e. $\mu_1,\mu_2\subset\partial B$).
\end{definition}
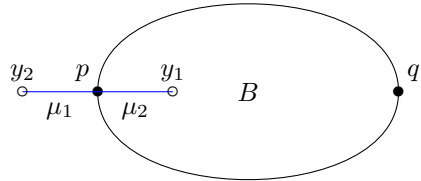
\begin{figure}[h]

\begin{center}
\begin{tikzpicture}
\node at (-2.2,.25) {$p$};
\node at (-2,0) {$\bullet$};
\node at (2.2,.25) {$q$};
\node at (2,0) {$\bullet$};
\node at (0,0) {$B$};
\draw (-2,0) to[out=90,in=90] (2,0);
\draw (-2,0) to[out=-90,in=-90] (2,0);

\draw[blue] (-3,0) to[out=0,in=180] (-1,0);
\node at (-3,0) {$\circ$};
\node at (-1,0) {$\circ$};
\node at (-3,0.25) {$y_2$};
\node at (-1,0.25) {$y_1$};
\node at (-2.5,-0.25) {$\mu_1$};
\node at (-1.5,-0.25) {$\mu_2$};

\end{tikzpicture}
\end{center}
\caption{A standard BM-configuration}
\end{figure}
Notice that, given a BM-configuration $(B,\mu,p)$ which is standard in the second sense, a very 
tiny isotopy of the bubble (which is allowed by Lemma \ref{injdevnbdofinjdevpath}) reduces $(B,\mu,p)$  
to 
a BM-configuration which is standard in the first sense. Namely in any projective coordinate we can 
push the developed image of the arc of bubbling slightly to the left or right of itself; when 
referring to a standard BM-configuration we will really always think of the first sense. We have 
the following characterisation.
\begin{lemma}\label{standardBMcharacterisation}
 Let $\sigma \in \mathcal{BM}_{(1,1),\rho}$ and let $(B,\mu,p)$ be a BM-configuration on it, such 
that $\sigma=Bub(\sigma_0,\beta)$ for some bubbleable arc $\beta \subset \sigma_0 \in 
\mathcal{M}_{0,\rho}$. Then $(B,\mu,p)$ is a standard 
BM-configuration if and only if $\mu$ induces an arc $\mu'$ on $\sigma_0$ such that the 
concatenation of $\beta$ and $\mu'$ is a bubbleable arc on $\sigma_0$. 
\end{lemma}
\begin{proof}
When we debubble $\sigma$ with respect to $B$ we naturally end up with the unbranched structure 
$\sigma_0$ endowed with a bubbleable arc $\beta$ such that  $Bub(\sigma_0,\beta)=\sigma$.
One of the two arcs contained in the embedded twin pair, let us say $\mu_2$, starts outside the 
bubble, hence its germ survives in $\sigma_0$, and we can try to analytically continue it to a path 
$\mu_0$ which has the same developed image of $\mu$. If the BM-configuration is standard then 
$\mu_2$ never meets the bubble, thus $\mu_0$ is a simple arc on $\sigma_0$, which does not meet 
$\beta$ away from $p$; in other words the concatenation of $\beta$ and $\mu_0$ is a simple arc on 
$\sigma_0$. Moreover the developed image of this arc is given by the concatenation of the developed 
image of $\partial B$ and $\mu$, which are disjoint. Thus this arc is bubbleable on $\sigma_0$. 
Conversely, if this arc is bubbleable, then when we perform the bubbling we can reconstruct the 
embedded twin pair $\mu$ by looking for the twin of $\mu_0$ inside the bubble. Since the whole 
$\beta\mu_0$ is bubbleable and we are bubbling only along the subarc $\beta$, we see that the 
developed image of the remaining part does not cross that of $\beta$. This means exactly that the 
twin starting inside the bubble will not leave it. Therefore the induced BM-configuration is 
standard. 
\end{proof}
The interest in standard BM-configurations is motivated by the following lemma.
\begin{lemma}\label{standardbm}
 Let $\sigma_0 \in \mathcal{BM}_{(1,1),\rho}$ and  $(B,\mu,p)$ be a standard BM-configuration on 
it; let $\sigma_t$ be the BPS obtained by moving branch points on $\sigma_0$ along $\mu$ up to 
time $t$, where $t\in [0,1]$ is a parameter along the developed image of $\mu$. Then $\sigma_t \in 
\mathcal{BM}_{(1,1),\rho}$ for all $t\in [0,1]$.
\end{lemma}
\begin{proof} 
This directly follows from the characterisation in Lemma \ref{standardBMcharacterisation} together with  
\cite[Lemma 2.9]{CDF}. In the above notations we have that 
$\sigma_t = $ $Move(\sigma_0,\mu^t)$  $ = $  $Bub(\sigma_0,\beta \mu'^t)$, where $\mu^t$ and $\mu'^t$ 
are the 
subarcs of $\mu$ and $\mu'$ respectively from time $0$ to time $t$. 
\end{proof}
We are now ready to prove the following result.
\begin{theorem}\label{bubopen}
 Let $\rho:\pi_1(S)\to\pslc$ be a non-elementary representation. Then $\mathcal{BM}_{(1,1),\rho}$ 
is open in $\mathcal{M}_{(1,1),\rho}$ (hence in $\mathcal{M}_{2,\rho}$).
\end{theorem}
\begin{proof} By Lemma \ref{standardbm} it is enough to show that given $\sigma_0 \in \mathcal{BM}_{(1,1),\rho}$ 
there is a small neighbourhood $U$ of it such that any structure in $U$ is obtained by moving 
branch 
points along a standard BM-configuration on $\sigma_0$. This easily follows from the fact the 
moving 
branch points gives a full neighbourhood of $\sigma_0$ in the moduli space (see Remark \ref{movingiscoord}), because local movement 
of branch points can always be performed along embedded twin pairs which are in standard 
BM-configuration with a given bubble on $\sigma_0$.
\end{proof}

Notice however that a priori more complicated BM-configurations might arise, which can not be used 
to move branch points preserving the bubble; namely if the embedded twin pair intersects the 
boundary of the bubble (or if this holds for their developed images), then moving branch 
points results in the break of the bubble: the aspiring bubbleable arc is either not embedded or 
not injectively developed. In this case it is not clear if it is possible to find another 
bubble on the spot.
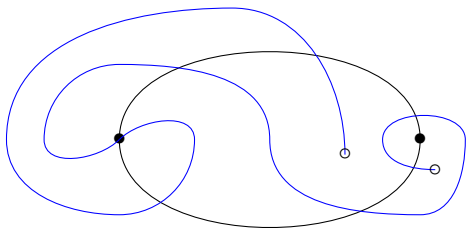
\begin{figure}[h]

\begin{center}
\begin{tikzpicture}
\node at (-2,0) {$\bullet$};
\node at (2,0) {$\bullet$};
\draw (-2,0) to[out=90,in=90] (2,0);
\draw (-2,0) to[out=-90,in=-90] (2,0);

\draw[blue] (-2,0) to[out=45,in=90] (-1,0) to[out=-90,in=0] (-2,-1) to[out=180,in=-90] (-3.5,0)
to[out=90,in=180] (-0.5,1.75) to[out=0,in=90] (1,-0.2);
  \node at (1,-0.2) {$\circ$};
 
\draw[blue] (-2,0) to[out=180+45,in=-90] (-3,0) to[out=90,in=180] (-2,1) to[out=0,in=90] (0,0)
to[out=-90,in=180] (2,-1) to[out=0,in=-90] (2.6,0) to[out=90,in=90] (1.5,0) to[out=-90,in=180] 
(2.2,-0.4) ;
  \node at (2.2,-0.4) {$\circ$};

\end{tikzpicture}
\end{center}
\caption{A non-standard BM-configuration}
\end{figure}

This heuristic argument can be made more precise by the following observation: moving branch points 
on a standard BM-configuration preserves the isotopy class (relative to endpoints) of the bubble; 
in particular it does not change the underlying unbranched structure. On the other hand it is not 
difficult to produce examples of movements of branch points which do not preserve the 
underlying unbranched structure. 
\begin{example}\label{ex_movhypbbltononsimplydev}
 Let $\rho:\pi_1(S)\to \pslc$ be a Fuchsian representation, $\beta$ be a bubbleable arc on the 
hyperbolic surface $\sigma_\rho=\hyp^2/\rho$ and $\sigma=Bub(\sigma_\rho,\beta)$. Notice that if 
$dev$ is a developing map for $\sigma$ and $x,y$ are its branch points, then it is not possible to 
find a pair of developed images $\widehat{x}$ and $\widehat{y}$ of them such that 
$\widehat{x}=\widehat{y}$; as we will say below, such a structure is said to be 
simply developed: the developed images of the branch point have disjoint $\rho$-orbits. On the 
other hand it is possible to move branch points on $\sigma$ along  suitable embedded twin pairs 
$\mu$ and $\nu$ with both endpoints inside the bubble in such a way that the resulting structure 
does not have this property (see Figure \ref{pic_movehyptononsimplydev}). This of course prevents 
the structure $Move(\sigma,\mu)$ from being a bubbling over $\sigma_\rho$. Of course the 
BM-configuration on $\sigma$ is not standard.
\begin{figure}[h]
\begin{center}
\begin{tikzpicture}[xscale=0.9]
 
\draw (0,0) to[out=90,in=180] (2,2) to[out=0,in=180] (4,1.5) to[out=0,in=180] (7,1.5) 
to[out=0,in=180] (9,2) to[out=0,in=90] (11,0) ;
\draw[yscale=-1] (0,0) to[out=90,in=180] (2,2) to[out=0,in=180] (4,1.5) to[out=0,in=180] (7,1.5) 
to[out=0,in=180] (9,2) to[out=0,in=90] (11,0) ;
\draw[xscale=0.8]  (1.6,0) to[out=65,in=180] (2.4,0.4) to[out=0,in=115] (3.2,0);
\draw[xscale=0.8]  (1.5,0.2) to[out=-75,in=180] (2.4,-0.3) to[out=0,in=-105] (3.3,0.2);
\draw[xshift=7cm,xscale=0.8]  (1.6,0) to[out=65,in=180] (2.4,0.4) to[out=0,in=115] (3.2,0);
\draw[xshift=7cm,xscale=0.8]  (1.5,0.2) to[out=-75,in=180] (2.4,-0.3) to[out=0,in=-105] (3.3,0.2);

\node at (3.5,0) {$\bullet$};
\node at (3,-0.25) {$x$};
\node at (7.5,0) {$\bullet$};
\node at (8,0.25) {$y$};
\draw[blue] (3.5,0) to[out=90,in=180] (5.5,1) to[out=0,in=90] (7.5,0);
\draw[blue,yscale=-1] (3.5,0) to[out=90,in=180] (5.5,1) to[out=0,in=90] (7.5,0);
\draw[blue] (4.5,0.25) to[out=90,in=180] (5.5,0.75) to[out=0,in=90] (6.5,0.25);
\draw[blue,yscale=-1] (4.5,0.25) to[out=90,in=180] (5.5,0.75) to[out=0,in=90] (6.5,0.25);

\draw[red] (3.5,0) to[out=0,in=120] (4.5,-1.5);
\draw[red,dashed] (4.5,-1.5) to[out=60,in=-120] (5,1.5);
\draw[red] (5,1.5) to[out=-30,in=120] (5.5,1);
\draw[red] (3.5,0) to[out=180,in=-90] (3,0.5) to[out=90,in=180] (4.25,0.85) to[out=0,in=170] 
(5.5,0.75) ;
\node at (5.5,1) {$*$};
\node at (5.5,0.75) {$*$};
\begin{scope}[yscale=-1,xscale=-1,xshift=-11cm]
 \draw[green] (3.5,0) to[out=0,in=120] (4.5,-1.5);
\draw[green,dashed] (4.5,-1.5) to[out=60,in=-120] (5,1.5);
\draw[green] (5,1.5) to[out=-30,in=120] (5.5,1);
\draw[green] (3.5,0) to[out=180,in=-90] (3,0.5) to[out=90,in=180] (4.25,0.85) to[out=0,in=170] 
(5.5,0.75) ;
\node at (5.5,1) {$*$};
\node at (5.5,0.75) {$*$};
\end{scope}

\foreach \x in {-1,0,1}
 \node at (12,\x) {$\bullet$};
\foreach \x in {-1,0,1}
 \node at (13,\x) {$*$};
\foreach \x in {-1,0,1}
 \node at (14,\x) {$\bullet$};
\foreach \x in {-1,0,1}
 \draw[blue] (12,\x) -- (14,\x);

 \node at (12,-.25) {$\widehat{x}$};
 \node at (14,.27) {$\widehat{y}$};

\draw[red] (12,0) to[out=45,in=180] (12.25,0.25) to[out=0,in=120] (13,-1);
\draw[green] (14,0) to[out=-45,in=0] (13.75,-0.25) to[out=180,in=-60] (13,1);

\node at (13,-2) {$\cp$};

\end{tikzpicture}

\end{center}
\caption{Figure for Example \ref{ex_movhypbbltononsimplydev}.}\label{pic_movehyptononsimplydev}
\end{figure}
\end{example}

\subsection{Taming developed images and avatars}
One of the main technical issues about $\cp$-structures is that the developing map is dramatically 
non-injective (already in the case of unbranched structures), hence it is quite difficult 
to control the relative behaviour of the developed images of some configuration of objects on 
the surface, even when the configuration is well behaved on the surface, as seen in Example 
\ref{ex_movhypbbltononsimplydev}.
\begin{definition}
Let $H\leq \pi_1(S)$ be a subgroup and $U\subset S$ be any subset. Let $\sigma$ be a 
BPS on $S$. We say that $U$ is \textbf{$H$-tame} (with respect to $\sigma$) if for some 
lift $\widetilde{U}$ of $U$ we have that a developing map for $\sigma$ is injective when 
restricted to $\cup_{h \in H} h. \widetilde{U}$. We will just say $U$ is tame if it is 
$\pi_1(S)$-tame.
\end{definition}
Notice that a tame simple arc is in particular bubbleable, and that a tame simple closed curve is 
in particular graftable as soon as the holonomy is loxodromic.
\begin{example}
 Any subset of a hyperbolic surface is tame, simply because the developing map is globally 
injective. More generally, if $\sigma$ is a \qf BPS and $C\subset \sigma^\pm$ is an unbranched 
geometric component, then any subset of $C$ is $\pi_1(C)$-tame, and any subset of the convex core 
of 
$C$ is tame.
\end{example}
Being able to control the collection of developed images of a given object on the surface (e.g. a 
curve) will not be enough in the following. For example, even if we start with a very well behaved 
structure $\sigma_0$ (e.g. a hyperbolic surface), when we perform a bubbling or a grafting we 
introduce in our structure $\sigma_0$ a  region $R$  whose full developed image is the whole 
model space $\cp$; as a result, inside $R$ we ``see'' a lot of developed images of any given subset 
$U\subset \sigma_0$. The following definition aims at making this more precise.
\begin{definition}
 Let $\rho:\pi_1(S)\to \pslc$ be a representation,  $\sigma \in \mathcal{M}_{k,\rho}$ 
and  $U\subset \sigma$ be any subset. An \textbf{avatar} of $U$ is any subset $V\subset \sigma$ 
such that there exist a lift $\widetilde{U}$ of $U$ and a lift $\widetilde{V}$ of $V$ such that 
$dev(\widetilde{U})=dev(\widetilde{V})$.
A structure $\sigma \in \mathcal{M}_{(1,1),\rho}$ is said to be \textbf{simply developed} if the 
two branch points are not avatars of each other.
\end{definition}

\begin{example}
If a structure has an injective developing map, then having the same developed image 
means being the same set, so that there are no non-trivial avatars. This happens for a hyperbolic 
surface, and more generally for the uniformizing structure $\sigma_\rho$ of a \qf representation 
$\rho:\pi_1(S)\to \pslc$.
\end{example}

In \qf holonomy we have a well-defined notion of size for subsets avoiding the real curve, 
which allows us to control the collection of avatars of a small set, as the following result shows. 
Let us denote by $sys(\rho)$ the systole of the uniformizing structure $\sigma_\rho$ or, 
equivalently, the minimum of the translation lengths of the elements in $\rho(\pi_1(S))$.

\begin{lemma}\label{smallsetistamedisjointavatar}
 Let $\rho:\pi_1(S)\to \pslc$ be Fuchsian and $\sigma \in \mathcal{M}_{0,\rho}$. Let 
$U\subset \sigma$ be a connected set with $diam(U)<sys(\rho)$ and which is $\pi_1$-trivial (i.e. 
$i_*(\pi_1(U))\subset \pi_1(S)$ is the trivial subgroup). Then $U$ sits inside a geometric 
component, it is tame and its avatars are disjoint.
\end{lemma}
\begin{proof}
Recall that when the holonomy is \qf there is a well defined hyperbolic metric on the complement of 
the real curve, which blows up in a neighbourhood of it; hence we can define a generalised 
path metric on the whole surface. Any connected subset of $\sigma$ which intersects the real curve 
must have infinite diameter with respect to this metric, because any path intersecting the real 
curve has infinite length. Therefore $U$ can not intersect the real curve, hence it is contained in 
some geometric component.\par
Since $U$ is $\pi_1$-trivial, it lifts homeomorphically to the universal cover. To prove tameness, 
assume that there are two lifts $\widetilde{U}_1$ and $\widetilde{U}_2$ which overlap once 
developed, i.e. $\exists \ x_i \in \widetilde{U}_i$ such that $dev(x_1)=dev(x_2)$. Let $\gamma \in 
\pi_1(S)$ be the unique deck transformation such that $\gamma \widetilde{U}_1=\widetilde{U}_2$. 
Then we have the following absurd chain of inequalities
$$sys(\rho)\leq d(\rho(\gamma) dev(x_1),dev(x_1))=d(\rho(\gamma) dev(x_1),dev(x_2))=$$
$$=d(dev(\gamma x_1),dev(x_2))\leq diam(dev (\widetilde{U}_2))=diam(U)<sys(\rho)$$
where $d$ denotes the hyperbolic distance on $\cp \setminus \rp$ and the last equality 
follows from 
the fact that the restriction of the developing map to each geometric component is an isometry.\par
Finally let us prove that the avatars in each geometric component are disjoint. Let $C$ be a 
geometric component, and choose a lift $\widetilde{C}$ of it and a lift $\widetilde{U}$ of $U$. The 
collection of avatars of $U$ in $C$ is given by
$$\pi_{|_{C}} \left( dev^{-1}_{|_{C}} \left( dev \left( \pi_1(S)\widetilde{U} \right) \right)  
\right) $$
So we want to prove that this is a disjoint collection. By tameness we know that the 
collection $ dev \left( \pi_1(S) \widetilde{U} \right)$ is disjoint, and the 
same is true for $dev^{-1}_{|_{C}} \left( dev \left( \pi_1(S)\widetilde{U} \right) \right)$, since 
the restriction of the developing map to each geometric component is a diffeomorphism and $ dev 
\left( \pi_1(S) \widetilde{U} \right)$ sits in the upper-half plane because $U$ is entirely 
contained in a positive geometric component. So we only 
need to prove that the projection $\pi$ does not overlap things too much. Let us introduce the 
following notation: if $\gamma \in \pi_1(S)$ then 
$$\gamma*\widetilde{U}:=dev^{-1}_{|_{C}} \left( dev \left( \gamma \widetilde{U} \right) \right)$$
With this notation what we want to prove now is that if there exist $\gamma_1,\gamma_2 \in 
\pi_1(S)$ such that $\pi\left( \gamma_1*\widetilde{U} \right) \cap \pi\left( \gamma_2*\widetilde{U} 
\right) \neq \varnothing$ then actually $\pi\left( \gamma_1*\widetilde{U} \right) = \pi\left( 
\gamma_2*\widetilde{U} \right)$. So let $x_i \in \gamma_i * \widetilde{U}$ such that 
$\pi(x_1)=\pi(x_2)$. Then $\exists \ \gamma \in \pi_1(C)$ such that $\gamma x_1=x_2$. If we develop 
these points we see that 
$$dev(x_2)=dev(\gamma x_1)=\rho(\gamma)dev(x_1)$$
and that $dev(x_2)\in dev(\gamma_2*\widetilde{U})=\rho(\gamma_2)dev(\widetilde{U})$ and 
$\rho(\gamma)dev(x_1) \in \rho(\gamma)dev(\gamma_1*\widetilde{U})=\rho(\gamma 
\gamma_1)dev(\widetilde{U})$. Since we already know that $U$ is tame, we can conclude that 
$\rho(\gamma \gamma_1)=\rho(\gamma_2)$, hence that $\gamma \gamma_1=\gamma_2$, because \qf 
representations are faithful. But then we have that
$$ \gamma_2*\widetilde{U}=(\gamma \gamma_1)*\widetilde{U}= dev^{-1}_{|_{C}} \left( dev \left( 
\gamma 
\gamma_1 \widetilde{U} \right) \right) =  dev^{-1}_{|_{C}} \left( \rho(\gamma) dev \left( \gamma_1 
\widetilde{U} \right) \right)$$
The last term is indeed equal to $\gamma dev^{-1}_{|_{C}} \left(  dev \left( \gamma_1 
\widetilde{U} \right) \right) $, because $\gamma \in \pi_1(C)$. So we have proved that 
$\gamma_2*\widetilde{U}=\gamma \left( \gamma_1*\widetilde{U} \right)$ for $\gamma\in \pi_1(C)$, 
which of course implies that $\pi\left( \gamma_1*\widetilde{U} \right) = \pi\left( 
\gamma_2*\widetilde{U} \right)$ as desired.
\end{proof}

Notice that the proof above  shows that in the collection $dev^{-1}_{|_{C}} ( dev 
( \pi_1(S)\widetilde{U} ) )$ either two elements differ by an automorphism of the 
universal cover $\pi:\widetilde{C}\to C$ and project to the same set on $C$, or they project to 
disjoint sets on $C$. In other words the avatars of $U$ in $C$ can be labelled by the cosets of 
$\pi_1(C)$ in $\pi_1(S)$; the index of $\pi_1(C)$ in $\pi_1(S)$ is  $1$ in the case of the 
uniformizing structure (where there are no non-trivial avatars, as already observed), and infinite 
otherwise, because in all the other cases any geometric component is a non-compact 
(incompressible) subsurface and free subgroups of surface groups have infinite index. We conclude 
with the following technical lemma which says that it is always possible to nicely isotope a 
bubbleable arc in order to minimise its intersections with a sufficiently small neighbourhood of 
its endpoints.
\begin{lemma}\label{isotopearcsoutsidesmallnbd}
 Let $\rho:\pi_1(S)\to \pslc$ be \qf and $\sigma \in \mathcal{M}_{0,\rho}$.
 Let $\beta \subset \sigma$ be a bubbleable arc with endpoints $x,y$ such that $x\not \in 
\sigma^\real$. Let $U\subset \sigma$ be a connected $\pi_1$-trivial  neighbourhood of $x$ 
with $diam(U)<sys(\rho)$  and not containing any avatar of $y$. Then there is an injectively 
developed isotopy (relative to $x$ and $y$) from $\beta$ to another bubbleable arc $\beta'$, such 
that $\beta'$ does not intersect any non-trivial avatar of $U$ and $\beta'\cap U$ is connected 
(i.e. $\beta'$ does not come back to $U$ after the first time it leaves it).
\end{lemma}
\begin{proof}
First of all notice that if $U$ does not contain avatars of $y$, then in particular $y$ is not an 
avatar of $x$. Moreover no avatar of $U$ contains avatars of $y$; in particular no avatar of $U$ 
contains $y$. 
We also know by Lemma \ref{smallsetistamedisjointavatar} that $U$ is geometric 
(i.e it avoids the real curve), tame and its avatars are disjoint. Since $U$ is geometric, for 
$\varepsilon >0$ small enough the $\varepsilon$-neighbourhood $\mathcal{N}_\varepsilon(U)$ of $U$ 
enjoys the same properties. \par
Let $\{U_i\}_{i\in I}$ be the collection of avatars of $U$ crossed by $\beta$. Going along $\beta$ 
from $x$ to $y$ we see that, apart from the initial segment starting at $x$ inside $U$, 
every time $\beta$ enters one of the $U_i$'s it crosses it and leaves it (this is exactly because 
no avatar contains the second endpoint $y$). Therefore we can  isotope all the arcs given by 
$\beta \cap U_i$ to arcs living in $\mathcal{N}_\varepsilon(U_i) \setminus U_i$, for each $i\in I$, 
without touching the first segment starting at $x$; since the chosen neighbourhood is tame this can 
be done in such a way that the isotopy is injectively developed. Since all the 
$\mathcal{N}_\varepsilon(U_i)$ are disjoint, this gives an 
isotopy on $\sigma$ from $\beta$ to an arc $\beta'$ which intersects the whole collection of 
avatars 
only in the initial segment starting at $x$ in $U$. It is still a bubbleable arc because it 
coincides with $\beta$ (which is bubbleable) outside the $\mathcal{N}_\varepsilon(U_i)$'s, and the 
deformations inside these sets do not produce any new intersection because 
$\mathcal{N}_\varepsilon(U)$ is tame.
\end{proof}
\begin{figure}[h]

\begin{center}
\begin{tikzpicture}

\fill[gray!20] (0,0) circle (0.5cm);
\fill[gray!5] (-2,0.5) circle (0.5cm);
\fill[gray!5] (1,-1) circle (0.5cm);

\node at (0,0) {$\bullet$};
\node at (-1,-1) {$\bullet$};
\draw (0,0) circle (0.5cm);
\draw (-2,0.5) circle (0.5cm);
\draw (1,-1) circle (0.5cm);

\node at (1,1) {\color{blue} $\beta$};
\draw[blue,thick] (0,0) to[out=0,in=-90] (2,0.15) to[out=90,in=0] (0,0.3) to[out=180,in=90] 
(-2,0.15)  to[out=-90,in=180] (0,-0.3) to[out=0,in=90] (1,-1)  to[out=-90,in=0] (0,-1.5)  
to[out=180,in=-90] (-1,-1) ;

\node (P) at (3,0) {};
\node (R) at (4,0) {};
\path[->,font=\scriptsize,>=angle 90]
(P) edge node[above]{} (R);

\fill[gray!20] (8,0) circle (0.5cm);
\fill[gray!5] (8-2,0.5) circle (0.5cm);
\fill[gray!5] (9,-1) circle (0.5cm);

\node at (8,0) {$\bullet$};
\node at (8-1,-1) {$\bullet$};
\draw (8,0) circle (0.5cm);
\draw (8-2,0.5) circle (0.5cm);
\draw (9,-1) circle (0.5cm);

\node at (9,1) {\color{blue} $\beta'$};

\draw[blue,thick,xshift=8cm] (0,0) to[out=0,in=-90] (2,0.15) to[out=90,in=0] (0.5,0.3)
to[out=110,in=70] (-0.5,0.3) 
to[out=180-20,in=-5] (-1.4,0.5) to[out=-90,in=0] (-2,-0.1)  to[out=-60,in=190] (-0.5,-0.3)   
to[out=-80,in=180+70] (0.5,-0.35)   to[out=-45,in=45+90] (0.7,-0.5)   to[out=180+30,in=180-35] 
(0.6,-0.5-0.9) to[out=180+20,in=-90] (-1,-1) ;

\end{tikzpicture}
\end{center}
\caption{Pushing an arc outside the avatars of a neighbourhood of one endpoint.}
\end{figure}
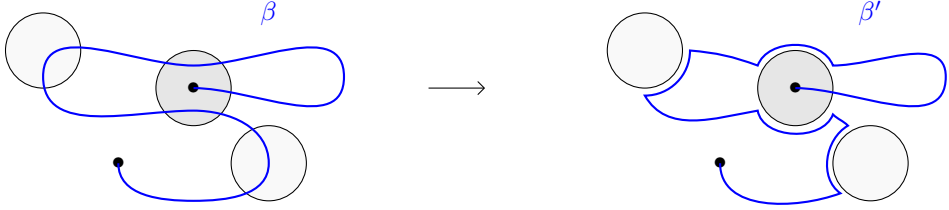
To get an intuition of what can go wrong, consider for instance the following picture; the 
second endpoint $y$ belongs to one of the avatars, hence there is no guarantee that the 
deformation that we want to perform is an injectively developed isotopy.
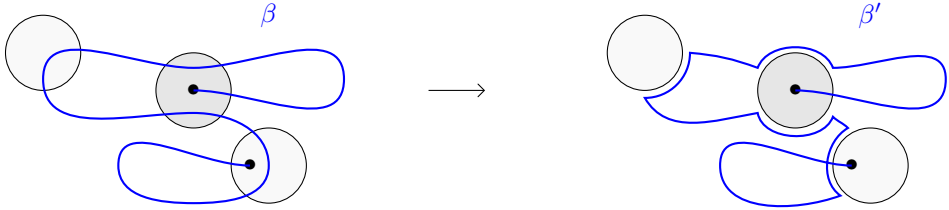
\begin{figure}[h]

\begin{center}
\begin{tikzpicture}

\fill[gray!20] (0,0) circle (0.5cm);
\fill[gray!5] (-2,0.5) circle (0.5cm);
\fill[gray!5] (1,-1) circle (0.5cm);

\node at (0,0) {$\bullet$};
\node at (0.75,-1) {$\bullet$};
\draw (0,0) circle (0.5cm);
\draw (-2,0.5) circle (0.5cm);
\draw (1,-1) circle (0.5cm);

\node at (1,1) {\color{blue} $\beta$};
\draw[blue,thick] (0,0) to[out=0,in=-90] (2,0.15) to[out=90,in=0] (0,0.3) to[out=180,in=90] 
(-2,0.15)  to[out=-90,in=180] (0,-0.3) to[out=0,in=90] (1,-1)  to[out=-90,in=0] (0,-1.5)  
to[out=180,in=-90] (-1,-1) to[out=90,in=180] (0.75,-1);

\node (P) at (3,0) {};
\node (R) at (4,0) {};
\path[->,font=\scriptsize,>=angle 90]
(P) edge node[above]{} (R);

\fill[gray!20] (8,0) circle (0.5cm);
\fill[gray!5] (8-2,0.5) circle (0.5cm);
\fill[gray!5] (9,-1) circle (0.5cm);

\node at (8,0) {$\bullet$};
\node at (8.75,-1) {$\bullet$};
\draw (8,0) circle (0.5cm);
\draw (8-2,0.5) circle (0.5cm);
\draw (9,-1) circle (0.5cm);

\node at (9,1) {\color{blue} $\beta'$};

\draw[blue,thick,xshift=8cm] (0,0) to[out=0,in=-90] (2,0.15) to[out=90,in=0] (0.5,0.3)
to[out=110,in=70] (-0.5,0.3) 
to[out=180-20,in=-5] (-1.4,0.5) to[out=-90,in=0] (-2,-0.1)  to[out=-60,in=190] (-0.5,-0.3)   
to[out=-80,in=180+70] (0.5,-0.35)   to[out=-45,in=45+90] (0.7,-0.5)   to[out=180+30,in=180-35] 
(0.6,-0.5-0.9) to[out=180+20,in=-90] (-1,-1)  to[out=90,in=180] (0.75,-1);;

\end{tikzpicture}
\end{center}
\caption{Avoiding avatars may result in self-intersections.}
\end{figure}

\subsection{Visible BM-configurations} 
This section is about a class of BM-con\-fig\-u\-ra\-tions with the property that, roughly 
speaking, the 
embedded twin pair survives after debubbling the structure, as in the proof of Lemma 
\ref{standardBMcharacterisation}; these should be thought as a strict 
generalisation of standard BM-configurations, which can still be dealt with by exploiting the 
underlying unbranched structure, where deformations are more easily defined and controlled. 

\begin{definition}
Let $\sigma_0 \in \mathcal{M}_{0,\rho}$, $\beta\subset \sigma_0$ a bubbleable arc and 
$\sigma=Bub(\sigma_0,\beta) \in \mathcal{M}_{2,\rho}$ with distinguished bubble $B$ coming from 
$\beta$. Let $p$ be a branch point of $\sigma$ and $\mu$ an embedded twin pair based at $p$ with 
developed image $\widehat{\mu}$. Notice that the germ of $\mu$ is well-defined on $\sigma_0$. We 
say that the BM-configuration $(B,\mu,p)$ is a \textbf{visible BM-configuration} if we can take the 
analytic continuation of this germ on $\sigma_0$ to obtain a properly embedded path $\mu_0$ on 
$\sigma_0$ which still develops to  $\widehat{\mu}$.
\end{definition}
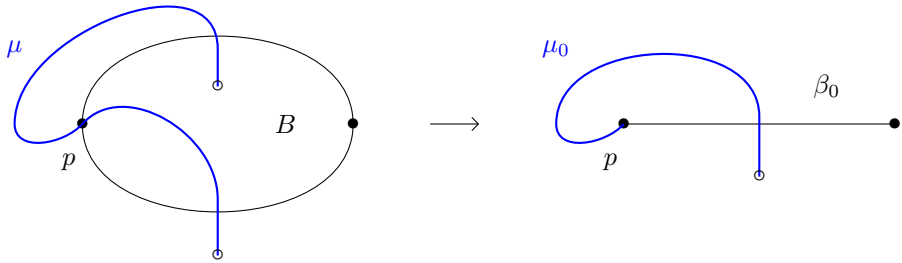
\begin{figure}[h]

\begin{center}
\begin{tikzpicture}[xscale=0.9]
\node at (-2,0) {$\bullet$};
\node at (-2.2,-0.5) {$p$};
\node at (-3,1) {\color{blue} $\mu$};
\node at (2,0) {$\bullet$};
\node at (1,0) {$B$};
\draw (-2,0) to[out=90,in=90] (2,0);
\draw (-2,0) to[out=-90,in=-90] (2,0);
 
\draw[blue,thick] (-2,0) to[out=180+45,in=-90] (-3,0) to[out=90,in=90] (0,1) to[out=-90,in=90] 
(0,0.5) ;
\node at (0,0.5) {$\circ$};

\draw[blue,thick] (-2,0) to[out=45,in=90] (0,-1) to[out=-90,in=90] (0,-1.75) ;
\node at (0,-1.75) {$\circ$};

\node (P) at (3,0) {};
\node (R) at (4,0) {};
\path[->,font=\scriptsize,>=angle 90]
(P) edge node[above]{} (R);

\draw (6,0) to[out=0,in=180] (10,0);
\node at (6,0) {$\bullet$};
\node at (10,0) {$\bullet$};
\node at (5.8,-0.5) {$p$};
\draw[blue,thick] (6,0) to[out=180+45,in=-90] (5,0) to[out=90,in=90] (8,0.1) to[out=-90,in=90] 
(8,-0.7) ;
\node at (8,-0.7) {$\circ$};
\node at (5,1) {\color{blue} $\mu_0$};
\node at (9,0.5) {$\beta_0$};

\end{tikzpicture}
\end{center}
\caption{A visible BM-configuration.}
\end{figure}
\begin{example}
 A standard BM-configuration is visible: as already observed in Lemma \ref{standardBMcharacterisation}, 
the boundary of the bubble and the embedded twin pair of a standard BM-configuration induce a pair 
of adjacent embedded arcs on the debubbled structure, whose concatenation is actually a bubbleable 
arc itself.
\end{example}
\begin{example}
 If $\sigma$ is a standard bubbling over a hyperbolic surface and we pick an embedded twin 
pair $\mu$ which intersects the real curve, then the resulting configuration is not visible: the 
debubbled structure is the uniformizing hyperbolic structure, which has no real curve, so there can 
be no path on it developing as needed; the analytic continuation of the germ of $\mu$ is an arc 
which wraps around the surface without converging to a compact embedded arc.
\end{example}

The next result shows that in \qf holonomy visible BM-configurations can be deformed to 
standard BM-configurations in a  controlled way.

\begin{proposition}\label{visibletostandard}
 Let $\rho:\pi_1(S)\to \pslc$ be Fuchsian, $\sigma_0 \in \mathcal{M}_{0,\rho}$, $\beta\subset 
\sigma_0$ a bubbleable arc. Let $x,y$ be the branch points of $\sigma = Bub(\sigma_0, \beta)$ and 
$B$ the bubble coming from $\beta$. Assume $\sigma$ is simply developed and $x\not \in 
\sigma^\real$. Let $K=\inf_{\gamma \in \pi_1(S)} d(dev(x),\rho(\gamma)dev(y))\in ]0,+\infty]$ and 
let $\mu$ be an embedded twin pair based at $x$ such that $(B,\mu,x)$ is a visible 
BM-configuration and the length of $\mu$ is less than $\min \{ sys(\rho),K\}$.
Then there is another bubble $B'\subset \sigma$ such that $Debub(\sigma,B')=\sigma_0$ and 
$(B',\mu,x)$ is a standard BM-configuration.
\end{proposition}
\begin{proof}
Since the BM-configuration is visible, after debubbling $\sigma$ we can define an arc $\mu_0$ on 
$\sigma_0$ starting at $x$ and developing as $\mu$. By hypothesis this arc is shorter than 
$sys(\rho)$ and $K$; in particular it can be put inside a connected contractible neighbourhood $U$ 
of $x$ with $diam(U)<sys(\rho)$ and which does not contain any avatar of $y$. By Lemma 
\ref{isotopearcsoutsidesmallnbd} there is an injectively developed isotopy from $\beta$ to a 
bubbleable arc $\beta'$ which avoids all non-trivial avatars of $U$ and intersects $U$ just once 
at the starting segment at $x$. Since this isotopy is injectively developed, bubbling $\sigma_0$ 
along $\beta'$ gives a structure isomorphic to $\sigma$ by Lemma \ref{bubblingonisotopicarcsareiso}.
Moreover the fact that $\mu \subset U$ and that $\beta'$ avoids all non-trivial avatars of $U$ and 
does not come back to it after the first time it leaves it implies that the concatenation of 
$\mu$ and $\beta'$ is a bubbleable arc; this is equivalent to saying that the resulting 
BM-configuration $(B',\mu,x)$ is standard by the characterisation in Lemma 
\ref{standardBMcharacterisation}.
\end{proof}
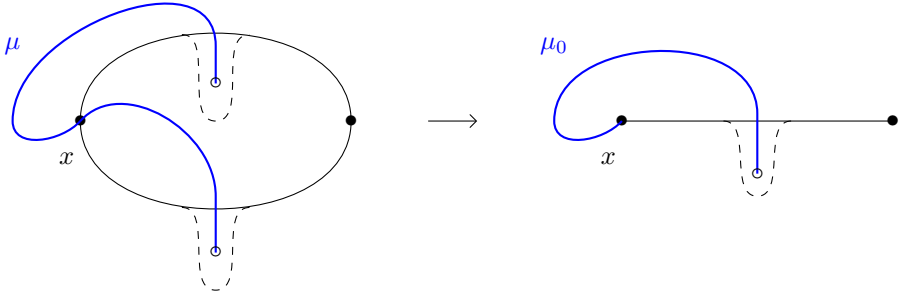
\begin{figure}[h]

\begin{center}
\begin{tikzpicture}[xscale=0.9]


\node at (-2,0) {$\bullet$};
\node at (-2.2,-0.5) {$x$};
\node at (-3,1) {\color{blue} $\mu$};
\node at (2,0) {$\bullet$};
\draw (-2,0) to[out=90,in=90] (2,0);
\draw (-2,0) to[out=-90,in=-90] (2,0);
 
\draw[blue,thick] (-2,0) to[out=180+45,in=-90] (-3,0) to[out=90,in=90] (0,1) to[out=-90,in=90] 
(0,0.5) ;
\node at (0,0.5) {$\circ$};

\draw[blue,thick] (-2,0) to[out=45,in=90] (0,-1) to[out=-90,in=90] (0,-1.75) ;
\node at (0,-1.75) {$\circ$};

\node (P) at (3,0) {};
\node (R) at (4,0) {};
\path[->,font=\scriptsize,>=angle 90]
(P) edge node[above]{} (R);

\draw (6,0) to[out=0,in=180] (10,0);
\node at (6,0) {$\bullet$};
\node at (10,0) {$\bullet$};
\node at (5.8,-0.5) {$x$};
\draw[blue,thick] (6,0) to[out=180+45,in=-90] (5,0) to[out=90,in=90] (8,0.1) to[out=-90,in=90] 
(8,-0.7) ;
\node at (8,-0.7) {$\circ$};
\node at (5,1) {\color{blue} $\mu_0$};

\draw[dashed] (7.5,0) to[out=0,in=180] (8,-1) to[out=0,in=180] (8.5,0);
\draw[dashed] (-.5,1.15) to[out=0,in=180] (0,0) to[out=0,in=180] (.5,1.15);
\draw[dashed] (-.5,-1.15) to[out=0,in=180] (0,-1.75-0.5) to[out=0,in=180] (.5,-1.15);

\end{tikzpicture}
\end{center}
\caption{Deforming a visible BM-configuration into a standard one.}
\end{figure}
\begin{rmk}

 The above result means that moving branch points on a given bubbling by a very small displacement 
(with respect to the representation) does preserve the [isotopy class of the] given 
bubble. In particular the underlying unbranched structure can be left unchanged throughout the 
movement. Notice that the hypothesis on the length is indeed necessary, as shown in Example 
\ref{ex_movhypbbltononsimplydev}: the BM-configuration therein is of course a visible one, but the 
embedded twin pairs are too long with respect both to the representation (the constant 
$sys(\rho)$ in Proposition \ref{visibletostandard}) and to the relative distance between the developed image 
of the branch points (the constant $K$ in Proposition \ref{visibletostandard}).
\end{rmk}
The condition of being visible is a bit obscure, if compared to that of being standard, in the 
sense that we have to debubble the structure to check visibility, and we do not have a simple 
characterisation as the one in Lemma \ref{standardBMcharacterisation} for standard BM-configurations; but 
visibility is always at least locally available at geometric branch points, as shown by the 
following result.
\begin{lemma}\label{smallvisible}
 Let $\rho:\pi_1(S)\to \pslc$ be Fuchsian, $\sigma_0 \in \mathcal{M}_{0,\rho}$, $\beta\subset 
\sigma_0$ a bubbleable arc such that $\sigma = Bub(\sigma_0, \beta) \in \mathcal{M}_{2,\rho}$ has 
a branch point $x$ not on the real curve. Let $\mu$ be an embedded twin pair based at $x$ of length 
smaller than $sys(\rho)$. Then the resulting BM-configuration is visible.
\end{lemma}
\begin{proof}
Let us fix a developed image $\widehat{x}$ of $x$ and $\widehat{\mu}$ for $\mu$. Since 
$l(\widehat{\mu})<sys(\rho)$, it is contained in a fundamental domain for 
$\rho$, and a fortiori in a fundamental domain for $\rho_{|_H}$, for any subgroup $H\leq \pi_1(S)$. 
Therefore it projects injectively to every quotient $\hyp^2/\rho(H)$.
Now consider the debubbled structure $\sigma_0$, and let $C$ be the geometric component of 
$\sigma_0$ containing $x$; since $\sigma_0$ is unbranched, $C$ is incompressible and carries a 
complete unbranched hyperbolic structure, so that the 
developing map induces an isometry $D_C:C\to \hyp^2/\rho(\pi_1(C))$, where 
$\widehat{\mu}$ projects injectively to an embedded arc. Pulling that arc back by $D_C$ gives the 
desired arc on $\sigma_0$, which proves the visibility of the BM-configuration.
\end{proof}
\begin{rmk}
 We want to remark that it is not possible to apply these ideas to a movement of a branch point 
which sits on the real curve. Indeed, here geometricity is used to produce neighbourhoods of the 
relevant objects which have disjoint avatars. On the other hand if a point belongs to the real 
curve, then any of its neighbourhoods will contain infinitely many avatars of both branch points, 
and actually of whatever object we want to consider. This follows from the fact that if $\Gamma$ is 
a Fuchsian group, then the collection of fixed points of its hyperbolic elements is dense in the 
limit set $\rp$.
\end{rmk}

\section{Bubbles everywhere}\label{s_bubbles}
In Theorem \ref{bubopen} we have proved that if $\rho:\pi_1(S)\to \pslc$ is a non-elementary 
representation, then the space of bubblings $\mathcal{BM}_{2,\rho}\subset 
\mathcal{M}_{2,\rho}$ is an open subspace of the moduli space $\mathcal{M}_{2,\rho}$. In this 
section we prove the main result of this paper, i.e. that in \qf holonomy it is also dense.
The strategy will be to consider a decomposition of the moduli space $\mathcal{M}_{2,\rho}$
obtained by looking at the position of the branch points with respect to the real curves: first of 
all we will show that if a piece of this decomposition contains a bubbling, then the bubblings fill 
a dense subspace in it, and then we will see that actually every piece of the decomposition 
contains a bubbling. The key step is the combination of Theorem \ref{k=2} above with the following result by Calsamiglia-Deroin-Francaviglia (see Theorem \cite[Theorem 7.1]{CDF}).
\begin{theorem}\label{t_k1cdf}
Let $S$ be a compact surface equipped with a BPS $\sigma$ having Fuchsian holonomy. Suppose there exist a positive and a negative component $C^+$ and $C^-$, with a common boundary component $l$ such that:
\begin{enumerate}
\item the index of $l$ is 1, and its holonomy is loxodromic;
\item the index of any other component of $\partial C^+$ and $\partial C^-$ vanishes;
\item $C^+$ and $C^-$ each contain a single simple branch point.
\end{enumerate}
Then, after possibly moving the branch points in the components $C^+$ and $C^-$, the branched projective structure   on $C^+\cup C^-$ is a bubbling.
\end{theorem}

Let us begin by defining the decomposition of $\mathcal{M}_{2,\rho}$ which we will use: recall from Definition
\ref{def_geomrealbranch} that if $\sigma$ is a BPS with \qf holonomy, then a branch point is 
geometric (respectively real) if it is outside (respectively on) the real curve.
Then $\sigma$ is said to be geometrically branched if it has no real branch points, while $\mathcal{M}_{2,\rho}^\real$ denotes the space of structures with at least one branch point.
\begin{lemma}
 The space $\mathcal{M}_{2,\rho}^\real$  has real codimension 1 in 
$\mathcal{M}_{2,\rho}$. 
\end{lemma}
\begin{proof}
This is a direct consequence of the fact that moving branch points on a structure provides local 
neighbourhoods for the manifold topology of $\mathcal{M}_{2,\rho}$ (see Remark \ref{movingiscoord}), and 
the fact that the real curve $\sigma^\real$ is an analytic curve on the surface.
\end{proof}
As a result the moduli space $\mathcal{M}_{2,\rho}$ decomposes as the union of the real 
hypersurfaces of $\mathcal{M}_{2,\rho}^\real$ and the remaining open pieces $\cup \mathcal{X}_i = 
\mathcal{M}_{2,\rho} \setminus \mathcal{M}_{2,\rho}^\real$ consisting of geometrically branched 
structures.
 \begin{definition}\label{def_realdecomp}
We will refer to the decomposition $\mathcal{M}_{2,\rho}=\mathcal{M}_{2,\rho}^\real \cup 
\bigcup_{i \in I}\mathcal{X}_{i}$ as the \textbf{real decomposition} of 
$\mathcal{M}_{2,\rho}$; any connected component $\mathcal{X}_i$ of $\mathcal{M}_{2,\rho}\setminus 
\mathcal{M}_{2,\rho}^\real$ will be called a geometric piece of the real decomposition of 
$\mathcal{M}_{2,\rho}$.
\end{definition}
Moving a geometrically branched structure $\sigma$ from a geometric piece $\mathcal{X}_i$ to an 
adjacent one $\mathcal{X}_j$ is a quite dramatic deformation, since it involves crossing 
$\mathcal{M}_{2,\rho}^\real$, i.e. moving branch points beyond the real curve $\sigma^\real$ of 
$\sigma$. This forces the combinatorial properties of the geometric decomposition to change 
abruptly; on the other hand, moving branch points on $\sigma$ inside their own geometric components 
keeps $\sigma$ inside the piece $\mathcal{X}_i$ it belongs to. This can actually be done in a 
quantitatively controlled way in order to preserve existing bubbles, as the following result shows.
\begin{theorem}\label{bubblinginconcpt}
 Let $\rho:\pi_1(S)\to \pslc$ be Fuchsian. Let $\mathcal{X} \subset \mathcal{M}_{2,\rho} 
\setminus \mathcal{M}_{2,\rho}^\real$ be a geometric piece of the real decomposition and let 
$\mathcal{Y}\subset \mathcal{X}$ be the subspace of simply developed structures. If $\mathcal{X}$ 
contains a bubbling, then every structure in $\mathcal{Y}$ is a bubbling; in particular 
$\mathcal{Y}=\mathcal{BY}$ is open and dense in $\mathcal{X}$.
\end{theorem}
\begin{proof}
First of all notice that $\mathcal{Y}$ is an open dense connected submanifold 
of $\mathcal{X}$, since its complement is a complex analytic subspace of complex (co)dimension 1.
As a consequence of Theorem \ref{bubopen}  $\mathcal{BX}$ is open; since $\mathcal{Y}$ is dense, it 
contains 
a bubbling too, i.e. $\mathcal{BY}$ is an open non-empty subset of $\mathcal{Y}$. We will prove 
that $\mathcal{BY}$ is also closed in $\mathcal{Y}$ and conclude by connectedness of 
$\mathcal{Y}$.\par
Let $\sigma_\infty \in \mathcal{Y}\cap \overline{\mathcal{BY}}$. By hypothesis the branch points 
$x_\infty$ and $y_\infty$ of $\sigma_\infty$ are outside the real curve of $\sigma_\infty$ and not 
avatar of each other. Fix any developed image $\widehat{x}_\infty$ of $x_\infty$ and 
$\widehat{y}_\infty$ of $y_\infty$; then define $K_\infty = \inf_{\gamma \in \pi_1(S)} 
d(\widehat{x}_\infty, \rho(\gamma)\widehat{y}_\infty)$. The distance here is the one induced by the 
hyperbolic metrics on the domain of discontinuity of $\rho$; $K_\infty$ is strictly positive since 
the branch points of $\sigma_\infty$ are not avatars, but can be $+\infty$ in the case they have 
opposite sign. Then let $A=\min\{sys(\rho),\frac{1}{3}K_\infty\} $.\par
Choose $L<A$ and consider the neighbourhood $\mathcal{N}_L(\sigma_\infty)$ of $\sigma_\infty$ in 
$\mathcal{Y}$ obtained by moving branch points by a distance $L < A$ (this is well defined since 
$\sigma_\infty$ is geometrically branched).
Since $\sigma_\infty$ is in the closure of $\mathcal{BY}$, $\mathcal{N}_L(\sigma_\infty)$ will 
contain a bubbling $\sigma \in \mathcal{BY}$. Let $\zeta$ be an 
embedded twin pair based at $x_\infty$ and $\xi$ be an embedded twin pair based at 
$y_\infty$ such that $Move(\sigma_\infty,\zeta,\xi)=\sigma$. By construction they can be chosen to 
have length smaller than $L$.
Let $\widehat{\zeta}$ and $\widehat{\xi}$ be the developed images based at 
$\widehat{x}_\infty,\widehat{y}_\infty$, and let $\sigma_0 \in \mathcal{M}_{0,\rho}$ and $\beta 
\subset \sigma_0$ be such that $\sigma=Bub(\sigma_0,\beta)$. Also let $x,y$ be the branch points 
of $\sigma$ corresponding to $x_\infty, y_\infty$ respectively.\par

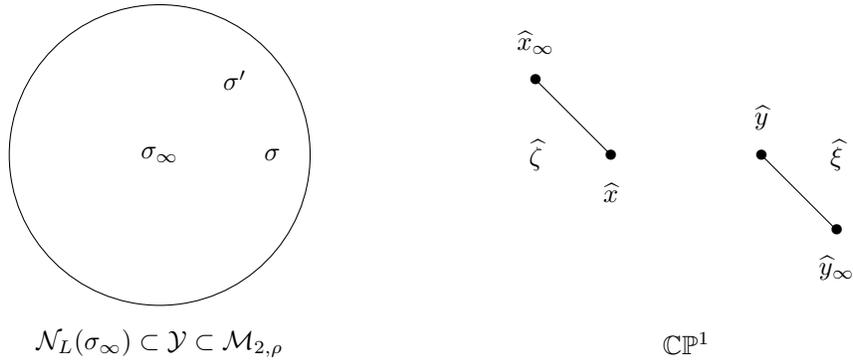
\begin{figure}[h]

\begin{center}
\begin{tikzpicture}


\draw (-3,0) circle (2cm);
\node at (-3,0) {$\sigma_\infty$};
\node at (-1.5,0) {$\sigma$};
\node at (-2,1) {$\sigma'$};
\node at (-3,-2.5) {$\mathcal{N}_L(\sigma_\infty)\subset \mathcal{Y}\subset \mathcal{M}_{2,\rho}$};

\node at (2,1) {$\bullet$}; \node at (2,1.5) {$\widehat{x}_\infty$};
\node at (3,0) {$\bullet$};\node at (3,-.5) {$\widehat{x}$};
\draw (2,1) -- (3,0);\node at (2,0) {$\widehat{\zeta}$};

\node at (5,0) {$\bullet$};\node at (5,.5) {$\widehat{y}$};
\node at (6,-1) {$\bullet$};\node at (6,-1.5) {$\widehat{y}_\infty$};
\draw (5,0) -- (6,-1);\node at (6,0) {$\widehat{\xi}$};

 \node at (4,-2.5) {$\cp$};

\end{tikzpicture}
\end{center}
 \caption{The neighbourhood $\mathcal{N}_L(\sigma_\infty)$ and the movement of points in $\cp$.}
\end{figure} 

We are now going to show that $\sigma_\infty$ is actually a bubbling over the same $\sigma_0$.
First of all notice that by Lemma \ref{smallvisible} both BM-configurations are visible, since both 
embedded twin pairs are shorter than the systole of the representation. Moreover, by definition of 
$A$, the two movements are independent from each other, i.e. commute; more precisely they do not 
interfere 
with each other in the sense that each twin pair avoids all avatars of the other twin pair. We 
begin by focusing at $x$; let us denote by $\widehat{x},\widehat{y}$ the developed images of $x,y$ 
which are seen at the endpoints of $\widehat{\zeta}$ and $\widehat{\xi}$. We have that for any 
$\gamma \in \pi_1(S)$
$$K_\infty \leq d( \widehat{x}_\infty,\rho(\gamma)\widehat{y}_\infty )\leq 
d(\widehat{x}_\infty,\widehat{x})+d(\widehat{x},\rho(\gamma)\widehat{y})+d(\rho(\gamma)\widehat{y},
\rho(\gamma)\widehat{y}_\infty)=2L+d(\widehat{x},\rho(\gamma)\widehat{y})$$
so that
$$d(\widehat{x},\rho(\gamma)\widehat{y})\geq K_\infty-2L$$
and we get 
$$\inf_{\gamma \in \pi_1(S)} (d(\widehat{x},\rho(\gamma)\widehat{y})  )\geq K_\infty-2L > 
3L-2L=L=l(\widehat{\zeta})$$
by definition of $A$. Therefore we can apply Proposition \ref{visibletostandard} and replace $\beta$ by a new 
bubbleable arc which is in standard BM-configuration on $\sigma$ with respect to $\zeta$. 
We now let $\sigma'=Move(\sigma,\zeta)$, which is still a bubbling over 
$\sigma_0$ by Lemma \ref{standardbm}.
We now want to use the same strategy again at $y$ to get back to $\sigma_\infty$; to do 
so, we just have to check that the movement is small enough with respect to the distance between 
the two branch points of $\sigma'$, which now develop to $\widehat{x}_\infty$ and $\widehat{y}$.
But this is easily checked: if $\gamma \in \pi_1(S)$ then
$$K_\infty \leq d( \widehat{x}_\infty,\rho(\gamma)\widehat{y}_\infty )\leq 
d(\widehat{x}_\infty,\rho(\gamma)\widehat{y})+d(\rho(\gamma)\widehat{y},
\rho(\gamma)\widehat{y}_\infty)=L+d(\widehat{x}_\infty,\rho(\gamma)\widehat{y})$$
so that
$$d(\widehat{x}_\infty,\rho(\gamma)\widehat{y})\geq K_\infty-L$$
and we get 
$$\inf_{\gamma \in \pi_1(S)} (d(\widehat{x}_\infty,\rho(\gamma)\widehat{y})  )\geq K_\infty-L > 
3L-L>L=l(\widehat{\xi})$$
So we can apply Proposition \ref{visibletostandard} again and replace the bubbleable arc with one which is in 
standard BM-configuration and safely move branch points along $\xi$. This movement results in our 
structure $\sigma_\infty$ and does not break the bubble by Lemma \ref{standardbm}. In other words this 
proves that $\sigma_\infty \in \mathcal{BY}$ (and indeed the underlying 
unbranched structure is the same as that of $\sigma$ and $\sigma'$), so that $\mathcal{BY}$ is 
closed.
\end{proof}
Let us denote by $k^\pm$ the number of positive and negative branch points of a structure as 
before. Notice that the value of $k^\pm$ is constant on every geometric piece of the real 
decomposition, so that it makes sense to say that a piece $\mathcal{X}$ has a given value of $k^+$. 
Combining all the results obtained so far, we can prove the following.
\begin{corollary}\label{bbleverywherek+=1}
 Let $\rho:\pi_1(S)\to \pslc$ be Fuchsian. Let $\mathcal{X}$ be a geometric 
piece of the real decomposition of $\mathcal{M}_{2,\rho}$ with $k^+=1$.
Then $\mathcal{X}=\mathcal{BX}$ i.e. is entirely made of bubblings.
\end{corollary}
\begin{proof}
First observe that every structure in $\mathcal{X}$ is simply developed, since branch 
points have different sign. Let $\sigma \in \mathcal{X}$. By Theorem \ref{k=2} it satisfies the hypothesis 
of \cite[Theorem 7.1]{CDF} (see Theorem \ref{t_k1cdf} above); therefore it is possible to move branch points inside their own 
geometric components so that a bubble appears, which proves that $\mathcal{X}$ contains a bubbling. 
The statement then follows from Theorem \ref{bubblinginconcpt}.
\end{proof}
We now have to care about geometric pieces of the real decomposition of $\mathcal{M}_{2,\rho}$ with 
$k^+=0,2$. We will prove that any such piece actually contains a bubbling. It should be said 
that the results in \cite[Proposition 8.1, Lemma 10.5-6]{CDF} imply that in some cases branch 
points can be moved inside their own geometric components so that a bubble appears, but it is not 
clear how to verify a priori when this occurs. Our strategy here will be to look for bubblings in 
the geometric pieces adjacent to $\mathcal{X}$ and drag them from there back into $\mathcal{X}$. In 
trying to do so, two problems occur: on one side if we naively take a bubbling in some piece 
adjacent to $\mathcal{X}$ and move branch points on it beyond the real curve, then it is quite 
difficult to control that we are actually moving to the chosen piece $\mathcal{X}$; on the other 
hand if we start with a structure  $\sigma \in \mathcal{X}$ and move branch points on it across the 
real curve to get to a bubbling, then it is quite difficult to check that when we move branch 
points 
to get back to $\sigma$ we do not break the bubble. Some lemmas are in order to guarantee that we 
can handle these issues. 
\begin{lemma}\label{k+1adjacentk+02}
Let $\mathcal{X}$ be a geometric piece of the real decomposition of $\mathcal{M}_{2,\rho}$ with 
$k^+=0$ or 2. Then there exists a geometric piece of the real decomposition $\mathcal{Y}$ adjacent 
to $\mathcal{X}$  and such that $k^+=1$.
\end{lemma}
\begin{proof}
This is just a reformulation of the results in \cite[\S 9]{CDF}, which say that it is always possible 
to move a branch point along a geodesic embedded twin pair crossing the real curve.
\end{proof}
We remark that in the process of moving a branch point towards the real curve with the 
techniques of \cite[\S 9]{CDF} a bubble might appear before actually crossing the real curve; this 
would be fine for us, since our ultimate goal now is to prove that $\mathcal{X}$ contains a 
bubbling; therefore we will forget about this detail in the following. 
The following lemma is needed to guarantee that it is always possible to go from one piece to an 
adjacent one by moving along a geodesic embedded twin pair.
\begin{lemma}\label{geodesictwinpair}
Let $\sigma \in \mathcal{M}_{2,\rho}$ be a geometrically branched BPS and $\mu=\{\mu_1,\mu_2\}$ 
an embedded twin pair on $\sigma$. Suppose that $\mu_i$ crosses 
$\sigma^\real$ at only one point $r_i$. Then there exists a geodesic embedded twin pair $\nu$ on 
$\sigma$ such that $Move(\sigma,\mu)=Move(\sigma, \nu)$.
\end{lemma}
\begin{proof}
Let $p$ be the base point of the embedded twin pair $\mu$ and $y_i$ be the endpoint of $\mu_i$.
By hypothesis the subarcs $\mu_i^1\subset \mu_i$ from $p$ to $r_i$ are entirely contained in a 
geometric component $C$. We let $\nu_i^1$ be the unique geodesic in $C$ from $p$ to $r_i$ which is 
isotopic to $\mu_i^1$ relatively to $\{p,r_i\}$. Then we can do the same in the adjacent components 
to obtain geodesic arcs $\nu_i^2$ isotopic to the subarcs $\mu_i^2\subset \mu_i$ from $r_i$ to 
$y_i$. 
\begin{figure}[h]

\begin{center}
\begin{tikzpicture}[scale=1,yscale=1.2]

\draw[xshift=10cm] (0,0) to[out=90,in=180] (2,1.5) to[out=0,in=180] (4,1) to[out=0,in=180] (7,1) 
to[out=0,in=180] (9,1.5) to[out=0,in=90] (11,0) ;
\draw[xshift=10cm,yscale=-1] (0,0) (0,0) to[out=90,in=180] (2,1.5) to[out=0,in=180] (4,1) 
to[out=0,in=180] (7,1) to[out=0,in=180] (9,1.5) to[out=0,in=90] (11,0) ;
\draw[xshift=10cm,xscale=0.8]  (1.6,0) to[out=65,in=180] (2.4,0.4) to[out=0,in=115] (3.2,0);
\draw[xshift=10cm,xscale=0.8]  (1.5,0.2) to[out=-75,in=180] (2.4,-0.3) to[out=0,in=-105] (3.3,0.2);
\draw[xshift=17cm,xscale=0.8]  (1.6,0) to[out=65,in=180] (2.4,0.4) to[out=0,in=115] (3.2,0);
\draw[xshift=17cm,xscale=0.8]  (1.5,0.2) to[out=-75,in=180] (2.4,-0.3) to[out=0,in=-105] (3.3,0.2);

\draw[xshift=10cm,red,thick] (4,1) to[out=180,in=90] (3.6,0) to[out=-90,in=180] (4,-1);
\draw[xshift=10cm,red,thick,dashed] (4,1) to[out=0,in=90] (4.4,0) to[out=-90,in=0] (4,-1);
\draw[xshift=13cm,red,thick] (4,1) to[out=180,in=90] (3.6,0) to[out=-90,in=180] (4,-1);
\draw[xshift=13cm,red,thick,dashed] (4,1) to[out=0,in=90] (4.4,0) to[out=-90,in=0] (4,-1);

\node at (15.5,0) {$\bullet$};
\node at (13.6,0) {$\bullet$};
\node at (12.8,0) {$\bullet$};
\draw[blue,thick] (15.5,0) .. controls (14.875,1)  .. (13.6,0);
\draw[blue,thick] (13.6,0) .. controls (13.2,-.5)  .. (12.8,0);
\draw[blue,dashed] (15.5,0) -- (12.8,0);

\node at (16.6,0) {$\bullet$};
\node at (18,0) {$\bullet$};
\draw[blue,dashed] (15.5,0) -- (18,0);
\draw[blue,thick] (15.5,0) .. controls (16,-.5)  .. (16.6,0);
\draw[blue,thick] (16.6,0) .. controls (17.3,1)  .. (18,0);

\node at (15.5,.5) {$p$};
\node at (14,-.5) {$r_1$};
\node at (17,-.5) {$r_2$};
\node at (12.8,.5) {$y_1$};
\node at (18,-.5) {$y_2$};

\end{tikzpicture}

\end{center}
 \caption{Straightening the embedded twin pair.}
\end{figure}
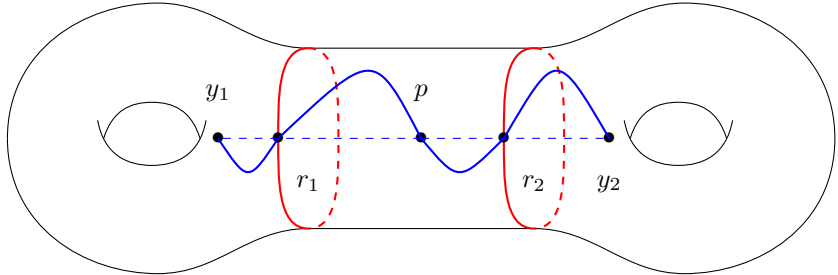
The concatenation of these paths gives rise to a pair of geodesic paths $\nu_i$ from $p$ 
to $y_i$ which are isotopic to $\mu_i$ relatively to $\{p,r_i,y_i\}$. Each geometric subarc 
$\mu_i^j$ 
is the geodesic representative of an embedded and injectively developed arc, hence it is 
embedded and injectively developed; moreover the two geometric subarcs of $\nu_i$ live in two 
adjacent component, hence their developed images are disjoint and $\nu$ is thus actually a geodesic 
embedded twin pair. The isotopy from $\mu$ to $\nu$ can be chosen to be an isotopy of embedded twin 
pairs, so that Lemma \ref{moveisotopic} applies and gives us that $Move(\sigma,\mu)=Move(\sigma, \nu)$.
\end{proof}

Notice that this result does not hold for paths which cross more components, because subarcs of 
$\nu$ 
contained in geometric components of the same sign can overlap once developed, even if $\mu$ is 
an embedded twin pair. We are now ready to prove that all the pieces of the real decomposition of 
$\mathcal{M}_{2,\rho}$ contain a bubbling. We will need the following terminology.
\begin{definition}
Let $\sigma \in \mathcal{M}_{k,\rho}$ be geometrically branched.
 Let $C\subset \sigma^\pm$ be a geometric component and $l\subset \partial C $ a real component in 
its boundary. We call the \textbf{peripheral geodesic} of $l$ in $C$ the unique geodesic 
representative $\gamma$ in the free homotopy class of $l$. The \textbf{end} of $l$ in $C$ is the 
connected component $E_l$ of $C \setminus \gamma$ which has $l$ in its boundary.
\end{definition}
It is shown in \cite[\S 3.3]{CDF} that ends are embedded open annuli, that ends associated to 
different real components are disjoint, and that if a geodesic enters an end, then it can not 
leave it and must necessarily reach the associated real curve.
\begin{theorem}\label{everycomponenthasbbl}
 Let $\rho:\pi_1(S)\to \pslc$ be a \qf representation. Let $\mathcal{X}$ be a geometric piece 
of the real decomposition of $\mathcal{M}_{2,\rho}$. Then $\mathcal{X}$ contains at least one 
structure which is a bubbling over an unbranched structure in $\mathcal{M}_{0,\rho}$.
\end{theorem}
\begin{proof}
If $\mathcal{X}$ has $k^+=1$ then this follows directly from Corollary \ref{bbleverywherek+=1}. So let us 
assume that $k^+=2$, the case $k^+=0$ being the same up to switching the signs of the branch 
points. We choose some $\sigma_1 \in \mathcal{X}$ and move branch points along an embedded twin 
pair $\mu$ to get to a structure $\sigma_2=Move(\sigma_1,\mu)$ in some adjacent piece $\mathcal{Y}$ 
with $k^+=1$, which can be done by Lemma \ref{k+1adjacentk+02}. By Theorem \ref{k=2} we know the combinatorial 
properties of the geometric decomposition of $\sigma_2$: all real curves are essential, one has 
index $1$ and the others have index $0$. Let us call $l$ the unique real curve of index $1$; the 
branch points $p^\pm$ live in the two geometric components $C^\pm$ adjacent to $l$. By construction 
we have an induced embedded twin pair $\nu$ at $p^-$ on $\sigma_2$ such that 
$Move(\sigma_2,\nu)=\sigma_1$. Here we know by \cite[Theorem 7.1]{CDF}  (see Theorem \ref{t_k1cdf} above) that we can move both branch 
points inside their own components to get to a structure $\sigma_3 \in \mathcal{Y}$ such that the 
peripheral geodesics of $l$ go through the  branch points $q^\pm$ of $\sigma_3$ with angles 
$\{\pi,3\pi\}$ and also such that it has a geodesic bubble $B$ (such a bubble can indeed be chosen 
in many ways, which will be exploited below). Of course we have an induced couple of embedded twin 
pairs $\zeta^\pm \subset \sigma_3$ based at $q^\pm$ such that 
$Move(\sigma_3,\zeta^+,\zeta^-)=\sigma_1$, and we would like to operate this movement of branch 
points on $\sigma_3$ without breaking the bubble $B$; unfortunately there is no reason why 
$(B,\zeta^\pm,q^\pm)$ should be a standard BM-configuration.\par

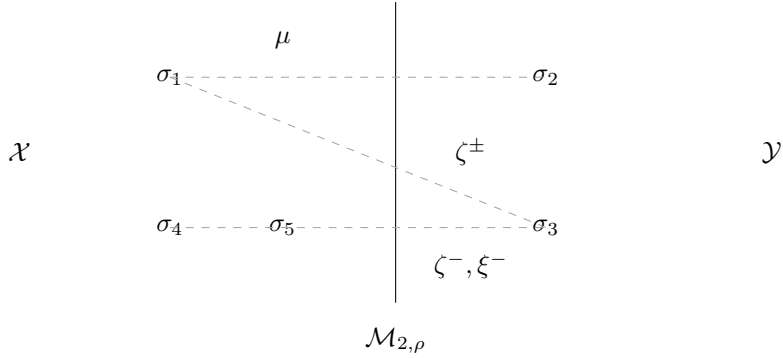
\begin{figure}[h]

\begin{center}
\begin{tikzpicture}


\draw (0,-2) -- (0,2) ;

\node at (-3,1) {$\sigma_1$};
\node at (-3,-1) {$\sigma_4$};
\node at (-1.5,-1) {$\sigma_5$};
\node at (2,1) {$\sigma_2$};
\node at (2,-1) {$\sigma_3$};
\node at (-5,0) {$\mathcal{X}$};
\node at (5,0) {$\mathcal{Y}$};
\node at (0,-2.5) {$\mathcal{M}_{2,\rho}$};

\draw[gray!70,dashed] (-3,1) -- (2,1) ;
\draw[gray!70,dashed] (-3,1) -- (2,-1) ;
\draw[gray!70,dashed] (-3,-1) -- (2,-1) ;

\node at (-1.5,1.5) {$\mu$};
\node at (1,0) {$\zeta^\pm$};
\node at (1,-1.5) {$\zeta^-,\xi^-$};

%
%
%

\end{tikzpicture}
\end{center}
 \caption{The structures $\sigma_1,\sigma_2,\sigma_3,\sigma_4$ and $\sigma_5$ involved in the 
proof.}
\end{figure}

However for our purposes we do not actually need to move branch points to go back to $\sigma_1$: 
it is enough to move to a structure in the same piece $\mathcal{X}$ without breaking the bubble 
$B$. Therefore we can forget about the embedded twin pair $\zeta^+$, since we only need to move 
$q^-$ to go back to that piece. Since $\zeta^-$ crosses the real curve just once, by Lemma
\ref{geodesictwinpair} we can replace it with a geodesic embedded twin pair $\xi^-$ which is such 
that $\sigma_4=Move(\sigma_3,\xi^-)=Move(\sigma_3,\zeta^-)\in \mathcal{X}$. As mentioned above, the 
bubble $B$ on $\sigma_3$ can be chosen in a quite free way, and our aim now is to prove that it 
is always possible to choose the bubble so that the BM-configuration $(B,\xi^-_{cut},q^-)$ is 
standard, for some suitable truncation $\xi^-_{cut}$ of the embedded twin pair $\xi^-$; of course 
we still have that $\sigma_5=Move(\sigma_3,\xi^-_{cut})\in \mathcal{X}$. \par
First of all we recall from \cite[\S 7]{CDF} that the real curve $l$ carries a natural action of the 
infinite cyclic group generated by $\rho(l)$ and a natural $\rho(l)$-invariant decomposition 
$l=\{0\}\cup l^+ \cup \{\infty\} \cup l^-$, corresponding to the decomposition of the limit set of 
$\rho$ given by the fixed points of $\rho(l)$; according to \cite[Proposition 7.8]{CDF} for any 
$u\in l^+$ we can find a geodesic bubble $B_u$ intersecting $l$ exactly at $u$ and $\rho(l)^{-1}u$.
Suppose we pick one of these geodesic bubbles $B_u$ and look at the situation on $C^-$, neglecting 
for a moment what happens beyond the real curve $l$. Since the embedded twin pair $\xi^-$ and the 
bubble $B_u$ are both geodesic, when one of the paths of $\xi^-$ enters the bubble it can never 
leave it, and must reach the real curve $l$. One of them, let 
us say $\xi^-_1$ starts inside $B_u$ (up to an arbitrarily small displacement of $u$), hence hits 
$l$ at some point $v_1$. If the BM-configuration $(B_u,\xi^-,q^-)$ is not already standard, it 
means that the twin $\xi^-_2$ starting outside $B_u$ goes somewhere around the surface and then 
comes back to intersect $B_u$ at some point $x$, and finally hits the real curve $l$ at some 
point $v_2$, distinct from $v_1$, because $\xi^-$ is an embedded twin pair.
Now, let us show that $v_2$ must live in $l^+$. 
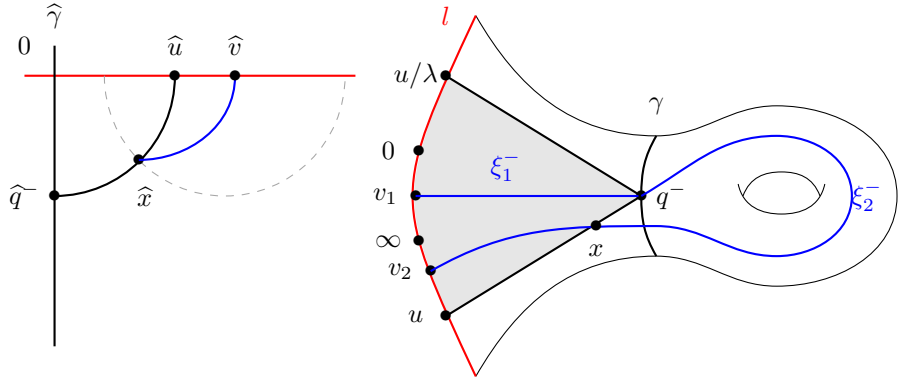
\begin{figure}[h]
\begin{center}
\begin{tikzpicture}[xscale=-1,scale=0.8]

\fill[gray!20] (4.25,0) -- (7.5,2) .. controls (8.2,0) .. (7.5,-2) -- (4.25,0);

\draw (0,0) to[out=90,in=180] (2,1.5) to[out=0,in=180] (4,1) to[out=0,in=240] (7,3);
\draw[yscale=-1] (0,0) to[out=90,in=180] (2,1.5) to[out=0,in=180] (4,1) to[out=0,in=240] (7,3);
\draw[xscale=0.8]  (1.6,0) to[out=65,in=180] (2.4,0.4) to[out=0,in=115] (3.2,0);
\draw[xscale=0.8]  (1.5,0.2) to[out=-75,in=180] (2.4,-0.3) to[out=0,in=-105] (3.3,0.2);

\draw[red,thick] (7,3) .. controls (8.4,0) .. (7,-3);
\draw[thick] (4,1) to[out=-60,in=90] (4.25,0) to[out=-90,in=60] (4,-1);

\node at (7.5,3) {\color{red} $l$};
\node at (4,1.5) { $\gamma$};
\node at (6.5,0.5) {\color{blue} $\xi_1^-$};
\node at (.5,0) {\color{blue} $\xi_2^-$};

\node at (4.25,0) {$\bullet$};
\node at (3.75,0) {$q^-$};
\node at (8,0) {$\bullet$};
\node at (8.5,0) {$v_1$};
\node at (7.5,2) {$\bullet$};\node at (8,2) {$u/\lambda$};
\node at (7.5,-2) {$\bullet$};\node at (8,-2) {$u$};

\draw[blue,thick] (4.25,0) -- (8,0);
\draw[thick] (4.25,0) -- (7.5,2);
\draw[thick] (4.25,0) -- (7.5,-2);

\node at (7.75,-1.25) {$\bullet$};\node at (8.25,-1.25) {$v_2$};
\node at (7.75+0.2,-0.75) {$\bullet$};\node at (8.25+0.2,-0.75) {$\infty$};
\node at (7.75+0.2,0.75) {$\bullet$};\node at (8.25+0.2,0.75) {$0$};

\draw[blue,thick] (4.25,0) to[out=150,in=0] (2,1) to[out=180,in=90] (.75,0) to[out=-90,in=180] 
(2,-1) to[out=0,in=180] (4,-0.5) to[out=0,in=150]  (7.75,-1.25);

\node at (5,-.5) {$\bullet$};\node at (5,-0.9) {$x$};


\draw[red,thick] (9,2) -- (14.5,2);
\draw[thick] (14,2.5) -- (14,-2.5);
\node at (14.5,2.5) {$0$};
\node at (14,3) {$\widehat{\gamma}$};
\node at (14,0) {$\bullet$};\node at (14.5,0) {$\widehat{q}^-$};
\node at (12,2) {$\bullet$};\node at (12,2.5) {$\widehat{u}$};
\draw[thick] (14,0) to[out=180,in=-90] (12,2);
\draw[dashed,gray!70] (14-0.828,2) to[out=-90,in=0] (14-2.828,0) to[out=180,in=-90] (14-4.828,2);

\node at (12.6,0.6) {$\bullet$};\node at (12.5,0) {$\widehat{x}$};
\draw[thick,blue] (12.6,0.6) to[out=180,in=-90] (11,2);
\node at (11,2) {$\bullet$};\node at (11,2.5) {$\widehat{v}$};

\end{tikzpicture}

\end{center}

 \caption{The configuration in $\cp$ and $C^-\subset \sigma_3$ when $B_u$ is the bubble orthogonal 
to the peripheral geodesic.}
\end{figure}
To do this, we choose $u$ so that the bubble $B_u$ 
is orthogonal at $q^-$ to the peripheral geodesic of $l$. Since $\xi^-_2$ is a geodesic from $q^-$ 
to $l$, once it enters the end relative to $l$ it constantly increases its distance from the 
peripheral geodesic; in particular, when it intersects the bubble at $x$ it forms an angle smaller 
than $\frac{\pi}{2}$  with the boundary of $B_u$. Since $u$ is in $l^+$, this forces $v_2\in l^+$ 
as well.
\begin{figure}[h]
\begin{center}
\begin{tikzpicture}[xscale=-1]

 \fill[gray!20] (4.25,0) -- (7.25,2.5) .. controls (8.2,0) .. (7.85,-1) -- (4.25,0);

\draw (0,0) to[out=90,in=180] (2,1.5) to[out=0,in=180] (4,1) to[out=0,in=240] (7,3);
\draw[yscale=-1] (0,0) to[out=90,in=180] (2,1.5) to[out=0,in=180] (4,1) to[out=0,in=240] (7,3);
\draw[xscale=0.8]  (1.6,0) to[out=65,in=180] (2.4,0.4) to[out=0,in=115] (3.2,0);
\draw[xscale=0.8]  (1.5,0.2) to[out=-75,in=180] (2.4,-0.3) to[out=0,in=-105] (3.3,0.2);

\draw[red,thick] (7,3) .. controls (8.4,0) .. (7,-3);
\draw[thick] (4,1) to[out=-60,in=90] (4.25,0) to[out=-90,in=60] (4,-1);

\node at (7.25,2.5) {$\bullet$};\node at (8,2.5) {$u'/\lambda$};
\node at (7.85,-1) {$\bullet$};\node at (8.25,-1.25) {$u'$};
\draw[thick] (4.25,0) -- (7.25,2.5);
\draw[thick] (4.25,0) -- (7.85,-1);


\node at (7.5,3) {\color{red} $l$};
\node at (4,1.5) { $\gamma$};
\node at (6.5,0.5) {\color{blue} $\xi_1^-$};
\node at (.5,0) {\color{blue} $\xi_2^-$};

\node at (4.25,0) {$\bullet$};
\node at (3.75,0) {$q^-$};
\node at (8,0) {$\bullet$};
\node at (8.5,0) {$v_1$};
\node at (7.75,-1.25) {$\bullet$};\node at (8,-1.8) {$v_2$};
\node at (7.75+0.2,-0.75) {$\bullet$};\node at (8.25+0.2,-0.75) {$\infty$};
\node at (7.75+0.2,0.75) {$\bullet$};\node at (8.25+0.2,0.75) {$0$};

\draw[blue,thick] (4.25,0) -- (8,0);
\draw[blue,thick] (4.25,0) to[out=150,in=0] (2,1) to[out=180,in=90] (.75,0) to[out=-90,in=180] 
(2,-1) to[out=0,in=180] (4,-0.5) to[out=0,in=150]  (7.75,-1.25);


%
%

\end{tikzpicture}

\end{center}

 \caption{A bubble in standard BM-configuration.}
\end{figure}
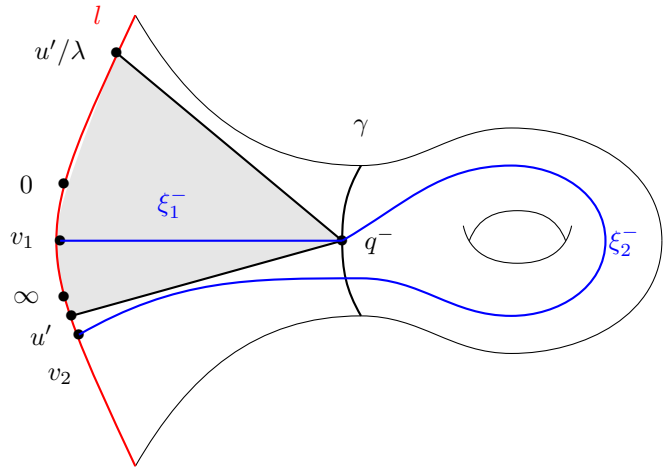
But then it is now possible to choose a different $u'$  in such a 
way that the arc $\alpha\subset l$ from $u'$ to $\rho(l)^{-1}u'$ containing $0$ and $\infty$ (i.e. 
the part of $l$ contained in $B_{u'}$) does not contain $v_2$. This choice guarantees that $v_2$ is 
outside the bubble $B_{u'}$, hence that $\xi^-_2$ does not intersect $B_{u'}$ before crossing the 
real curve $l$. We have no tools to control what happens beyond $l$, but we can truncate 
$\xi^-$ to a sub-embedded twin pair $\xi^-_{cut}$ which ends beyond $l$ and which is in 
standard BM-configuration with respect to the bubble $B_{u'}$. By Lemma \ref{standardbm}
$\sigma_5=Move(\sigma_3,\xi^-_{cut})$ is still a bubbling. But we can clearly keep moving branch 
points on $\sigma_5$ along what is left of $\xi^-$ to reach the structure 
$\sigma_4=Move(\sigma_3,\xi^-)$, which, as we already know, lives in the same piece 
$\mathcal{X}$ containing $\sigma_1$. Since this movement does not cross the real curve, the 
structure $\sigma_5$ lives in $\mathcal{X}$ too, which proves that $\mathcal{X}$ contains a 
bubbling.
\end{proof}

We can finally prove the main result. 
\begin{theorem}\label{mainbubblingthm}
Let $\rho:\pi_1(S)\to \pslc$ be a Fuchsian representation. Then any simply developed structure with at most one real branch point is a bubbling. In particular the space of bubblings is a connected, open and dense subspace of full measure in $\mathcal{M}_{2,\rho}$.
\end{theorem}
\begin{proof}
At first let $\sigma$ be a geometrically branched and simply developed structure. Since its branch 
points are outside the real curve it belongs to some geometric piece $\mathcal{X}$ of the real 
decomposition. By Theorem \ref{everycomponenthasbbl} we know that $\mathcal{X}$ contains a bubbling. 
Moreover $\sigma$ avoids the subspace of $\mathcal{X}$ made of non-simply developed structures. 
Then 
by Theorem \ref{bubblinginconcpt} $\sigma$ is a bubbling. In the case $\sigma$ has one real branch point, 
we can perform a movement of that branch point to go from $\sigma$ to some structure $\sigma'$ in 
some geometric piece of the real decomposition with $k^+=1$. Then the previous arguments apply 
verbatim, because the isotopy in Lemma \ref{geodesictwinpair} fixes the points of intersection between 
the 
embedded twin pair and the real curve, so that we are able to pick a bubble on $\sigma'$ and move 
back to $\sigma$ as in Theorem \ref{everycomponenthasbbl}. The subspace of structures left outside by this 
approach is the union of the  subspaces of non-simply developed structures and the one of 
structures with both branch points on the real curve; each of them has real 
codimension 2 in $\mathcal{M}_{2,\rho}$, which is a connected manifold of real dimension 4 by 
\cite{CDF} (see for instance Theorem \ref{t_complex_manifold} above), so that the last statement follows.
\end{proof}

\subsection{Walking around the moduli space with bubblings}
As a consequence of the results obtained in this paper we get a generically positive answer in our 
setting to the question asked by Gallo-Kapovich-Marden as Problem 12.1.2 in \cite{GKM}, i.e. if any 
two BPS with the same holonomy are related by a sequence of grafting, degrafting, bubbling and 
debubbling. 
More precisely Theorem \ref{mainbubblingthm} shows that, if $\sigma$ and $\tau$ are a generic 
pair of BPS with at most two branch points and a fixed \qf holonomy, then we can apply one 
debubbling to each of them (if needed), to reduce to a pair of unbranched structures 
$ \sigma_0$ and $\tau_0$ with the same holonomy. By Goldman's theorem in \cite{GO} we can then 
apply $m$ degraftings on $\sigma_0$ to obtain the uniformizing structure $\sigma_\rho$ and then $n$ 
graftings on $\sigma_\rho$ to obtain $\tau_0$, for suitable $m,n\in\nat$.\par
\begin{figure}[h]

\begin{center}
\begin{tikzpicture}
\node at (-6,0) {$\sigma$};

\draw[->] (-5.5,0) to[out=0,in=180] (-3.5,0);
\node at (-4.5,.5) {$1$ debub};

\node at (-3,0) {$\sigma_0$};

\draw[->] (-2.5,0) to[out=0,in=180] (-.5,0);
\node at (-1.5,.5) {$m$ degraft};

\node at (0,0) {$\sigma_\rho$};

\draw[->] (.5,0) to[out=0,in=180] (2.5,0);
\node at (1.5,.5) {$n$ graft};

\node at (3,0) {$\tau_0$};

\draw[->] (3.5,0) to[out=0,in=180] (5.5,0);
\node at (4.5,.5) {$1$ bub};

\node at (6,0) {$\tau$};

\end{tikzpicture}
\end{center}
\end{figure}
Actually it is possible to do even better, since we can remove the need for degraftings; by the 
proof of \cite[Theorem 11]{CDF2}, there exists a simple closed geodesic  $\gamma$ on $\sigma_\rho$ 
such that $\sigma_\gamma=Gr(\sigma_\rho,\gamma)$ can be obtained by $m'$ graftings on $\sigma_0$ and
$\tau_0$ can be obtained by $n'$ graftings on $\sigma_\gamma$, for suitable $m',n'\in\nat$.

Finally, according to \cite[Theorem 5.1]{CDF} every simple grafting can be realised by a sequence 
of one bubbling and one debubbling. This implies the following, which shows that it is generically 
possible to 
move around the moduli space only via bubblings and debubblings.
\begin{corollary}\label{onlybubanddebub}
Let $\rho:\pi_1(S)\to \pslc$ be Fuchsian.
 There is a connected, open and dense subspace  $\mathcal{B}\subset \mathcal{M}_{2,\rho}$ such 
that if 
$\sigma,\tau \in  \mathcal{M}_{0,\rho} \cup \mathcal{B}$ then $\tau$ is obtained from $\sigma$ by a 
finite sequence of bubblings and debubblings.
\end{corollary}
\begin{figure}[h]

\begin{center}
\begin{tikzpicture}
\node at (-6,0) {$\sigma$};

\draw[->] (-5.5,0) to[out=0,in=180] (-3.5,0);
\node at (-4.5,.5) {$1$ debub};

\node at (-3,0) {$\sigma_0$};

\draw[->] (-2.5,0) to[out=0,in=180] (-.5,0);
\node at (-1.5,.5) {$m'$ bub};
\node at (-1.5,-.5) {$m'$ debub};

\node at (0,0) {$\sigma_\gamma$};

\draw[->] (.5,0) to[out=0,in=180] (2.5,0);
\node at (1.5,.5) {$n'$ bub};
\node at (1.5,-.5) {$n'$ debub};

\node at (3,0) {$\tau_0$};

\draw[->] (3.5,0) to[out=0,in=180] (5.5,0);
\node at (4.5,.5) {$1$ bub};

\node at (6,0) {$\tau$};

\end{tikzpicture}
\end{center}
\end{figure}
Notice that the length of this sequence depends on the choice of the unbranched structures 
$\sigma_0$ and $\tau_0$ (i.e. the choice of the bubbles on $\sigma$ and $\tau$), which are not 
uniquely determined: a BPS with two branch points can in general be realised as a bubbling over 
different unbranched structures along different arcs. This phenomenon is outside the point of view 
of this paper, which was concerned with the preservation of the underlying unbranched structure 
during all the deformations, and is dealt with in a separate paper by the author (see 
\cite{R}).

\printbibliography

\end{document}